\documentclass[12pt]{article}
\usepackage{amsfonts}

\def\hybrid{\topmargin 0pt      \oddsidemargin 0pt
        \headheight 0pt \headsep 0pt
        \voffset=-0.5cm
        \textwidth 6.5in        
        \textheight 9in         
        \marginparwidth 0.0in
        \parskip 5pt plus 1pt   \jot = 1.5ex}
\catcode`\@=11
\def\marginnote#1{}

\newcount\hour
\newcount\minute
\newtoks\amorpm
\hour=\time\divide\hour by60
\minute=\time{\multiply\hour by60 \global\advance\minute by-\hour}
\edef\standardtime{{\ifnum\hour<12 \global\amorpm={am}%
        \else\global\amorpm={pm}\advance\hour by-12 \fi
        \ifnum\hour=0 \hour=12 \fi
        \number\hour:\ifnum\minute<10 0\fi\number\minute\the\amorpm}}
\edef\militarytime{\number\hour:\ifnum\minute<10 0\fi\number\minute}

\def\draftlabel#1{{\@bsphack\if@filesw {\let\thepage\relax
   \xdef\@gtempa{\write\@auxout{\string
      \newlabel{#1}{{\@currentlabel}{\thepage}}}}}\@gtempa
   \if@nobreak \ifvmode\nobreak\fi\fi\fi\@esphack}
        \gdef\@eqnlabel{#1}}
\def\@eqnlabel{}
\def\@vacuum{}
\def\draftmarginnote#1{\marginpar{\raggedright\scriptsize\tt#1}}
\def\draftlabel#1{{\@bsphack\if@filesw {\let\thepage\relax
   \xdef\@gtempa{\write\@auxout{\string
      \newlabel{#1}{{\@currentlabel}{\thepage}}}}}\@gtempa
   \if@nobreak \ifvmode\nobreak\fi\fi\fi\@esphack}
        \gdef\@eqnlabel{#1}}
\def\@eqnlabel{}
\def\@vacuum{}
\def\draftmarginnote#1{\marginpar{\raggedright\scriptsize\tt#1}}

\def\draft{\oddsidemargin -.5truein
        \def\@oddfoot{\sl preliminary draft \hfil
        \rm\thepage\hfil\sl\today\quad\militarytime}
        \let\@evenfoot\@oddfoot \overfullrule 3pt
        \let\label=\draftlabel
        \let\marginnote=\draftmarginnote
   \def\@eqnnum{(\theequation)\rlap{\kern\marginparsep\tt\@eqnlabel}%
\global\let\@eqnlabel\@vacuum}  }

\def\numberbysection{\@addtoreset{equation}{section}
        \def\theequation{\thesection.\arabic{equation}}}

\def\underline#1{\relax\ifmmode\@@underline#1\else
        $\@@underline{\hbox{#1}}$\relax\fi}

\def\titlepage{\@restonecolfalse\if@twocolumn\@restonecoltrue\onecolumn
     \else \newpage \fi \thispagestyle{empty}\c@page\z@
        \def\thefootnote{\fnsymbol{footnote}} }

\def\endtitlepage{\if@restonecol\twocolumn \else  \fi
        \def\thefootnote{\arabic{footnote}}
        \setcounter{footnote}{0}}  
\relax

\numberbysection
\hybrid

\def\beq{\begin{equation}}
\def\eeq{\end{equation}}
\def\bea{\begin{eqnarray}}
\def\eea{\end{eqnarray}}
\def\p{\partial}
\def\G{\Gamma}
\def\g{\gamma}

\def\z{\zeta}
\def\L{{\cal L}}

\def\C{{\cal C}}

\def\a{\alpha}
\def\b{\beta}
\def\e{\varepsilon}
\def\l{\lambda}
\def\f{\varphi}

\def\A{{\cal A}}

\def\D{{\cal D}}
\def\F{{\cal F}}

\def\L{{\cal L}}

\def\O{{\cal O}}

\def\U{{U_0}}

\def\dim{{\rm dim}}
\def\res{{\rm res}}

\def\TC{{\cal T}^t}
\def\TD{{\cal T}^{\nu}}
\def\wt{\widetilde}

\def \matrix #1 {\left(\begin{array}{cc} #1 \end{array}\right)}
\newtheorem{theo}{Theorem}[section]

\newtheorem{cor}{Corollary}[section]
\newtheorem{lem}{Lemma}[section]

\begin{document}
\begin{titlepage}
\title{Characterizing Jacobians via trisecants of the Kummer Variety}
\author{I.Krichever \thanks{Columbia University, New York, USA and
Landau Institute for Theoretical Physics, Moscow, Russia; e-mail:
krichev@math.columbia.edu. Research is supported in part by National Science
Foundation under the grant DMS-04-05519.}}

\date{May 22, 2006\footnote{the revised version February 28, 2007}}

\maketitle

\begin{abstract} We prove Welters' trisecant conjecture:
an indecomposable principally polarized abelian variety $X$ is the Jacobian of a curve
if and only if there exists a trisecant of its Kummer variety $K(X)$.
\end{abstract}

\end{titlepage}
\section{Introduction}

Welters' remarkable trisecant conjecture formulated first in \cite{wel1}
was motivated by Gunning's celebrated theorem (\cite{gun1})
and by another famous conjecture: the Jacobians of curves are exactly
the indecomposable principally polarized abelian varieties whose theta-functions
provide explicit solutions of the so-called KP equation. The latter
was proposed earlier by Novikov and was unsettled at the time of the Welters' work.
It was proved later by T.Shiota \cite{shiota} and until recently has remained
the most effective solution of the classical Riemann-Schottky problem.

Let $B$ be an indecomposable symmetric matrix with positive definite imaginary part.
It defines an indecomposable principally polarized abelian variety
$X=\mathbb C^g/\Lambda$, where  the lattice $\Lambda$ is generated by the basis vectors
$e_m\in \mathbb C^g$ and the column-vectors $B_m$ of $B$.
The Riemann theta-function $\theta(z)=\theta(z|B)$ corresponding to $B$
is given by the formula
\beq\label{teta1}
\theta(z)=\sum_{m\in \mathbb{Z}^g} e^{2\pi i(z,m)+\pi i(Bm,m)},\ \
(z,m)=m_1z_1+\cdots+m_gz_g\, .
\eeq
The Kummer variety $K(X)$ is an image of the Kummer map
\beq\label{kum}
K:Z\in X\longmapsto
\{\Theta[\e_1,0](Z):\cdots:\Theta[\e_{2^g},0](Z)\}\in \mathbb{CP}^{2^g-1}
\eeq
where $\Theta[\e,0](z)=\theta[\e,0](2z|2B)$ are level two theta-functions
with half-integer characteristics $\e$.

A trisecant of the Kummer variety is a projective line which meets
$K(X)$ at least at three points. Fay's well-known trisecant formula
\cite{fay} implies that if $B$ is a matrix of $b$-periods of normalized holomorphic
differentials on a smooth genus $g$ algebraic curve $\G$, then
a set of three arbitrary distinct points $A_1,A_2,A_3$ on $\G$
defines a {\it one-parametric
family} of trisecants parameterized by a fourth point of the curve $A_4\neq A_1,A_2,A_3$.
In \cite{gun1} Gunning proved under certain non-degeneracy assumptions
that the existence of such a family of trisecants
characterizes Jacobian varieties among indecomposable principally polarized
abelian varieties.

Gunning's geometric characterization of the Jacobian
locus was extended by Welters who proved that the Jacobian locus can
be characterized by the existence of a formal one-parametric family of flexes of the
Kummer varieties \cite{wel1,wel}. A flex of the Kummer variety
is a projective line which is tangent to $K(X)$ at some point up to order 2.
It is a limiting case of trisecants when the three intersection points come together.

In \cite{arb-decon} Arbarello and De Concini showed that the Welters' characterization
is equivalent to an infinite system of partial differential equations  representing
the so-called KP hierarchy, and proved that only a finite number of
these equations is sufficient.
In fact, the KP theory and the earlier results of Burchnall, Chaundy and
the author \cite{ch1,ch2,kr1,kr2} imply that the Jacobian locus is characterized
by the first $N=g+1$ equations of the KP hierarchy, only. Novikov's conjecture
that just the first equation ($N=1$!) of the hierarchy is sufficient for
the characterization of the Jacobians is much stronger.
It is equivalent to the statement that the Jacobians are characterized by the existence
of length $3$ formal jet of flexes.

In \cite{wel1} Welters formulated the question: {\it if the Kummer-Wirtinger variety $K(X)$ has {\it one}
trisecant, does it follow that $X$ is a Jacobian ?} In fact, there are
three particular cases of the Welters' conjecture, which are independent and have to be
considered separately. They correspond to three possible configurations of
the intersection points
$(a,b,c)$ of $K(X)$ and the trisecant:

(i) all three points coincide $(a=b=c)$,

(ii) two of them coincide $(a=b\neq c)$;

(iii) all three intersection points are distinct
$(a\neq b\neq c\neq a)$.

The affirmative answer to the the first particular
case (i) of the Welters' question was obtained in the author's previous work
\cite{kr-schot}. (Under various additional assumptions in various  forms
it was proved earlier in \cite{mar,kr3,flex}).
The aim of this paper is to prove, using the approach proposed in \cite{kr-schot},
the two remaining cases of the trisecant conjecture.
It seems that the approach is very robust and can be applied to the
variety of Riemann-Schottky-type problems. For example, in \cite{prym} it was
used for the characterization of principally polarized Prym
varieties of branched  covers.

Our first main result is the following statement.
\begin{theo} An indecomposable, principally polarized abelian variety $(X,\theta)$
is the Jacobian of a smooth curve of genus g if and only if
there exist non-zero $g$-dimensional vectors
$U\neq A \, (\bmod\,  \Lambda),\, V$, such that one of the following equivalent conditions
holds:

$(A)$  The differential-difference equation
\beq\label{laxd}
\left(\p_t-T+u(x,t)\right)\psi(x,t)=0, \ \ T=e^{\p_x}
\eeq
is satisfied for
\beq\label{u}
u=(T-1)v(x,t),\ \ v=-\p_t\ln\theta (xU+tV+Z)
\eeq
and
\beq\label{p}
\psi={\theta(A+xU+tV+Z)\over \theta(xU+tV+Z)}\, e^{xp+tE},
\eeq
where $p,E$ are constants and $Z$ is arbitrary.

\medskip
$(B)$ The equations
\beq\label{gr1}
\p_{V}\Theta[\e,0]\left((A-U)/2\right)-e^{p}\Theta[\e,0]\left((A+U)/2\right)
+E\Theta[\e,0]\left((A-U)/2\right)=0,
\eeq
are satisfied for all $\e\in {1\over 2}Z_2^g$. Here and below $\p_V$ is the constant
vector field on $\mathbb{C}^g$ corresponding to the vector $V$.

\medskip
$(C)$ The equation
\beq\label{cm7}
\p_V\left[\theta(Z+U)\,\theta(Z-U)\right]\p_V\theta(Z)=
\left[\theta(Z+U)\,\theta(Z-U)\right]\p^2_{VV}\theta(Z)\ ({\rm mod}\, \theta)
\eeq
is valid on the theta-divisor $\Theta=\{Z\in X\,\mid\, \theta(Z)=0\}$.
\end{theo}
Equation (\ref{laxd}) is one of the two auxiliary linear problems for the $2D$ Toda lattice
equation, which can be regarded  as a discretization of the KP equation. The idea to use
it for the characterization of Jacobians was motivated
by (\cite{kr-schot}), and the author's earlier work with Zabrodin (\cite{zab}), where
a connection of the theory of elliptic solutions of the $2D$ Toda lattice equations
and the theory of the elliptic Ruijsenaars-Schneider system was established.
In fact, Theorem 1.1 in a slightly different form was proved in (\cite{zab})
under the additional assumption that the vector $U$ {\it spans an elliptic curve} in $X$.

\medskip

The equivalence of $(A)$ and $(B)$ is a direct corollary of the addition formula for
the theta-function. The statement $(B)$ is the second particular case of
the trisecant conjecture: the line in $\mathbb{CP}^{2^g-1}$ passing through the points
$K((A-U)/2)$ and $K((A+U)/2)$ of the Kummer variety is tangent
to $K(X)$ at the point $K((A-U)/2)$.

The "only if" part of $(A)$ follows from the author's construction
of solutions of the $2D$ Toda lattice equations \cite{kr-toda}.
The statement $(C)$ is actually what we use for the proof of the theorem.
It is stronger than $(A)$. The implication $(A)\to (C)$ does not require
the explicit theta-functional formula for $\psi$. It is enough to require only that
equation (\ref{laxd}) with $u$ as in (\ref{u}) has {\it local meromorphic} in $x$
solutions which are holomorphic outside  the divisor $\theta(Ux+Vt+Z)=0$.

To put it more precisely, let $\tau(x,t)$ be a holomorphic function of $x$ in some domain $\D$,
where it has a simple root $\eta(t)$. If $\tau(\eta(t)\pm 1,t)\neq 0$, then the condition
that equation (\ref{laxd}) with $u=(T-1)v$, where $v(x,t))=-\p_t \ln\tau(x,t)$,
has a meromorphic solution with the only pole in $\D$ at $\eta$ implies
\beq\label{cm5}
\ddot \eta=\dot \eta\left[2v_0(t)-v(\eta+1,t)-v(\eta-1,t)\right]\,,
\eeq
where ``dots" stands for the $t$-derivatives and $v_0$ is the coefficient of the Laurent
expansion of $v(x,t)$ at $\eta$, i.e.
\beq\label{e1}
v(x,t)={\dot \eta\over x-\eta}+v_0(t)+O((x-\eta)).
\eeq
Formally, if we represent $\tau$ as an infinite product,
\beq\label{roots}
\tau(x,t)=c(t)\prod_i(x-x_i(t)),
\eeq
then equation (\ref{cm5}) can be written as the infinite system of equations
\beq\label{cm6}
\ddot x_i=\sum_{j\neq i} \dot x_i\dot x_j\left[{2\over (x_i-x_j)}-{1\over (x_i-x_j+1)}
-{1\over (x_i-x_j-1)}\right]\,.
\eeq
If $\tau$ is a rational, trigonometric or elliptic polynomial, then
the system (\ref{cm6}) coincides with the equations of
motion for the rational, trigonometrical or elliptic Ruijsenaars-Schneider systems,
respectively. Equations (\ref{cm6}) are analogues of the equations derived in \cite{flex}
and called in \cite{kr-schot} the formal Calogero-Moser system.

Simple expansion of $\theta$ at the points of its divisor $Z\in\Theta: \theta(Z)=0$
shows that for $\tau=\theta(Ux+Vt+Z)$ equation (\ref{cm5})
is equivalent to (\ref{cm7}).

The proof of the theorem goes along the same lines as the proof of Theorem 1.1 in
\cite{kr-schot}. In order to stress the similarity we almost literally copy some
parts of \cite{kr-schot}.
At the beginning of the next section we derive equations (\ref{cm5})
and show that they are sufficient conditions for the {\it local}
existence of {\it formal wave solutions}.
The formal wave solution of equation (\ref{laxd}) is a solution of the form
\beq\label{ps}
\psi(x,t,k)=k^xe^{kt}\left(1+\sum_{s=1}^{\infty}\xi_s(x,t)\,k^{-s}\right)\,.
\eeq
The ultimate goal is to show the existence of the wave solutions such that
coefficients of the series (\ref{ps}) have the form $\xi_s=\xi_s(Ux+Vt+Z)$, where
\beq\label{int}
\xi_s(Z)={\tau_s(Z)\over\theta(Z)}\,,
\eeq
and $\tau_s(Z)$ is a holomorphic function.
The functions $\xi_s$ are defined recursively by differential-difference equations
$(T_U-1)\,\xi_{s+1}=\p_V\xi_s+u\xi_s$, where $T_U=e^{\p_U}$ and $\p_U$ is a constant
vector-field defined by the vector $U$.

In the case of differential equations the  cohomological arguments, that
are due to Lee-Oda-Yukie, can be
applied for an attempt to glue local solutions into the global ones
(see details in \cite{shiota, arbarello}). These
arguments were used in \cite{shiota} and revealed that the core
of the problem in the proof of Novikov's conjecture is
a priori nontrivial cohomological obstruction for the global solvability
of the corresponding equations. The hardest part of the Shiota's work was
the proof that the certain bad locus $\Sigma\subset \Theta$, which controls
the obstruction, is empty \footnote{The author is grateful to Enrico Arbarello
for an explanation of these deep ideas and a crucial role of the singular locus $\Sigma$,
which helped him to focus on the heart of the problem.}.

In the difference case there is no analog of the cohomological arguments and we use
a different approach. Instead of {\it proving} the global existence
of solutions we, to some extend, {\it construct} them by defining
first their residue on the theta-divisor. It turns out that the residue
is regular on $\Theta$ outside the {\it singular locus} $\Sigma$
which is the maximal $T_U$-invariant subset of $\Theta$, i.e.
$\Sigma=\bigcap_{k\in\mathbb Z}T_U^k\Theta$.

As in \cite{kr-schot}, we don't prove directly that the bad locus is empty.
Our first step is to construct certain wave solutions outside the bad locus.
We call them $\l$-periodic wave solutions.
They are defined uniquely up to $t$-independent $T_U$-invariant factor.
Then we show that for each $Z\notin \Sigma$ the $\l$-periodic wave solution
is a common eigenfunction of a commutative ring $\A^Z$ of
ordinary difference operators.
The coefficients of these operators are independent of
ambiguities in the construction of $\psi$. For the generic $Z$ the ring $\A^Z$ is maximal
and the corresponding spectral curve $\G$ is $Z$-independent. The correspondence
$j:Z\longmapsto \A^Z$ and the results of the works \cite{mum,kr-dif}, where
a theory of rank 1 commutative rings of difference operators was developed,
allows us to make the next
crucial step and prove the global existence of the wave function.
Namely, on $(X\setminus\Sigma)$ the wave function can be globally defined as
the preimage $j^*\psi_{BA}$ under $j$ of the Baker-Akhiezer function on $\G$ and
then can be extended on $X$ by usual Hartogs' arguments. The global existence of the
wave function implies that $X$ contains an orbit of the
KP hierarchy, as an abelian subvariety. The orbit is isomorphic to the generalized
Jacobian $J(\G)={\rm Pic}^0(\G)$ of the spectral curve (\cite{shiota}).
Therefore, the generalized Jacobian
is compact. The compactness of $J(\G)$ implies that the spectral curve is smooth
and the correspondence $j$ extends by linearity and defines the isomorphism
$j: X\to J(\G)$.

In the last section we  present the proof  of the last "fully discrete" case of the
trisecant conjecture.

\begin{theo} An indecomposable, principally polarized abelian variety $(X,\theta)$
is the Jacobian of a smooth curve of genus g if and only if
there exist non-zero $g$-dimensional vectors
$U\neq V\neq A \neq U\, (\bmod \Lambda)$ such that one of the following equivalent
conditions holds:

$(A)$  The difference equation
\beq\label{laxdd}
\psi(m,n+1)=\psi(m+1,n)+u(m,n)\psi(m,n)
\eeq
is satisfied for
\beq\label{ud}
u(m,n)={\theta((m+1)U+(n+1)V+Z)\,\theta(mU+nV+Z)\over
\theta(mU+(n+1)V+Z)\,\theta((m+1)U+nV+Z)}
\eeq
and
\beq\label{pd}
\psi(m,n)={\theta(A+mU+nV+Z)\over \theta(mU+nV+Z)}\, e^{mp+nE},
\eeq
where $p,E$ are constants and $Z$ is arbitrary.

\medskip
$(B)$ The equations
\beq\label{gr1d}
\Theta[\e,0]\left({A-U-V\over 2}\right)+e^{p}\Theta[\e,0]\left({A+U-V\over 2}\right)
=e^E\Theta[\e,0]\left({A+V-U\over 2}\right),
\eeq
are satisfied for all $\e\in {1\over 2}Z_2^g$.

\medskip
$(C)$ The equation
\beq\label{cm7d}
\theta(Z+U)\,\theta(Z-V)\,\,\theta(Z-U+V)+\theta(Z-U)\,\theta(Z+V)\,\,\theta(Z+U-V)=0\ ({\rm mod}\, \theta)
\eeq
is valid on the theta-divisor $\Theta=\{Z\in X\,\mid\, \theta(Z)=0\}$.
\end{theo}
Equation (\ref{laxdd}) is one of the two auxiliary linear problems for the so-called
bilinear discrete Hirota equation  (BDHE). The "only if" part of $(A)$ follows
from the author's work \cite{krdd}.
Under the  assumption that the vector $U$ {\it spans an elliptic curve} in $X$,
theorem 1.2 was proved in \cite{bete}, where the connection of the
elliptic solutions of BDHE and the so-called elliptic nested Bete ansatz equations
was established.

\section{$\l$-periodic  wave solutions}

To begin with, let us show that equations (\ref{cm5}) are necessary for the
existence of a meromorphic solution of equation (\ref{laxd}), which is holomorphic
outside of the theta-divisor.

Let $\tau(x,t)$ be a holomorphic function of the variable $x$ in some {\it
translational invariant} domain $\D=T\D\subset \mathbb C$, where $T:x\to x+1$.
We assume that $\tau$ is a smooth function of the parameter $t$.  Suppose that
$\tau$ in $\D$ has a simple root $\eta(t)$ such that
\beq\label{xi}
\tau(\eta(t)+1,t)\,\tau(\eta(t)-1,t)\neq 0.
\eeq
\begin{lem}
If equation (\ref{laxd}) with the potential
$u=(T-1)v$, where $v=-\p_t\ln\tau(x,t)$ has a meromorphic in $\D$ solution $\psi(x,t)$,
with the simple pole at $x=\eta$, and regular at $\eta-1$, then equation (\ref{cm5}) holds.
\end{lem}
{\it Proof.}
Consider the Laurent expansions of $\psi$ and $v$ in the neighborhood of one of $\eta$:
\beq\label{ue}
v={\dot \eta\over x-\eta}+v_0+\ldots\,,
\eeq
\beq\label{psie}
\psi={\a\over x-\eta}+\b+\ldots\,.
\eeq
(All coefficients in these expansions are smooth functions of the variable $t$).
Substitution of (\ref{ue},\ref{psie}) in (\ref{laxd}) gives an infinite system of
equations. We use only the following three of them.

The vanishing of the residue at $\eta$ of the left hand side of (\ref{laxd}) implies
\beq\label{eq1}
\dot\a=\dot \eta \b+\a(v_0-v(\eta+1,t)).
\eeq
The vanishing of the residue and the constant terms of the Laurent expansion
of (\ref{laxd}) at $\eta-1$ are equivalent to the equations
\beq\label{eq2}
\a=\dot \eta\psi(\eta-1,t)\,,
\eeq
\beq\label{eq3}
\p_t \psi(\eta-1,t)=\b+[v(\eta-1,t)-v_0]\psi(\eta-1,t)\,.
\eeq
Taking the $t$-derivative of (\ref{eq2}) and using equations (\ref{eq1}, \ref{eq3})
we get (\ref{cm5}).

Let us show that equations (\ref{cm5}) are sufficient for the existence of {\it local}
meromorphic wave solutions which are holomorphic outside of the zeros of $\tau$.
In the difference case a notion of local solutions needs some clarification.

In the lemma below we assume that a translational invariant domain $\D$ is a disconnected
union of small discs, i.e.
\beq\label{dd}
\D=\cup_{i\in \mathbb Z}\  T^iD_0,\ \ D_0=\{x\in \mathbb C\, |\,x-x_0|<1/2\}.
\eeq
\begin{lem} Suppose that $\tau(x,t)$ is holomorphic in a domain $\D$ of the form
(\ref{dd}) where it has simple zeros, for which condition (\ref{xi}) and
equation (\ref{cm5}) hold.
Then there exist meromorphic wave solutions of equation (\ref{laxd}) that have
simple poles at zeros of $\tau$ and are holomorphic everywhere else.
\end{lem}
{\it Proof.} Substitution of (\ref{ps}) into (\ref{laxd}) gives a recurrent system of
equations
\beq\label{xis}
(T-1)\xi_{s+1}=\dot\xi_s+u\xi_s.
\eeq
Under the assumption that $\D$ is a disconnected union of small disks, $\xi_{s+1}$
can be defined as an arbitrary meromorphic function in $D_0$ and then
extended on $\D$ with the help of (\ref{xis}). If $\eta$ is a zero of $\tau$, then
in this way we get a meromorphic function $\xi_{s+1}$, which a priori has poles at the
points $\eta_k=\eta-k$ for all non-negative $k$. Our goal is to prove by induction
that in fact (\ref{xis}) has meromorphic solutions with simple poles only at the zeros
of $\tau$.

Suppose that $\xi_s$ has a simple pole at $x=\eta$
\beq\label{5}
\xi_s={r_s\over x-\eta}+r_{s0}+r_{s1}(x-\eta)+\cdots\,.
\eeq
The condition that $\xi_{s+1}$ has no pole at $\eta-1$ is equivalent to the equation
\beq\label{res}
R_s=\dot \eta\xi_s(\eta-1,t),
\eeq
where by definition
\beq\label{resa}
R_s=r_s(v(\eta+1,t)-v_0)+\dot \eta r_{s0}-\dot r_s\,.
\eeq
Equation (\ref{res}) with $R_s$ given by (\ref{resa}) is our induction assumption.
We need to show that the next equation holds also.
From (\ref{xis}) it follows that
\beq\label{6}
r_{s+1}=\dot \eta\xi_s(\eta-1,t)\,,
\eeq
\beq\label{7}
\xi_{s+1}(\eta-1,t)-r_{s+1,0}=-\dot \xi_s(\eta-1,t)(v(\eta-1,t)-v_0)\xi_s(\eta-1,t)\,.
\eeq
These equations imply
\beq\label{8}
R_{s+1}=\dot \eta\xi_{s+1}(\eta-1,t)-
(\ddot \eta+\dot \eta(v(\eta+1,t)+v(\eta-1,t)-2v_0))\xi_s(\eta-1,t)
\eeq
and the lemma is proved.

If $\xi_{s+1}^0$ is a particular solution of (\ref{xis}), then the general solution
is of the form $\xi_{s+1}(x,t)=c_{s+1}(x,t)+\xi_{s+1}^0(x,t)$
where $c_{s+1}(x,t)$ is $T$-invariant function of the variable
$x$, and an {\it arbitrary} function of the variable $t$.
Our next goal is to fix a {\it translation-invariant} normalization of
$\xi_s$.

Let us show that in the periodic case $v(x+N,t)=v(x,t)=-\p_t\tau(x,t), \, N\in \mathbb Z$,
the periodicity condition for $\xi_{s+1}(x+N,t)=\xi_{s+1}(x,t)$
uniquely defines $t$-dependence of the functions $c_s(x,t)$ (compare with
the normalization of the Bloch solutions of differential equations used in
\cite{kp}).
Assume that $\xi_{s-1}$ is known  and satisfies the condition that there exists a periodic
solution $\xi_s^0$ of the corresponding equation.
Let $\xi_{s+1}^*$ be a solution of (\ref{xis})
for $\xi_s^0$. Then the function $\xi_{s+1}^0=\xi_{s+1}^*+x\p_tc_s+c_sv$
is a solution of (\ref{xis}) for $\xi_s=\xi_s^0+c_s$.
A choice of $T$-invariant function $c_s(x,t)$ does not affect
the periodicity property of $\xi_s$, but it does affect the periodicity in $x$ of
the function $\xi_{s+1}^0(x,t)$.
In order to make  $\xi_{s+1}^0(x,t)$ periodic,
the function $c_s(x,t)$ should satisfy the linear differential equation
\beq\label{kp4}
N\p_t c_s(x,t)+\xi_{s+1}^*(x+N,t)-\xi_{s+1}^*(x,t)=0.
\eeq
That defines $c_s(x,t)$ uniquely up to a $t$-independent $T$-invariant function
of the variable $x$.

\bigskip
In the general case, when $U$ is not a point of finite order in $X$, the solution of the
normalization problem for the coefficients of the wave solutions requires the global
existence of these coefficients along certain affine subspaces.

Let $Y_U$ be the Zariski closure
of the group $\{Uk \mid k\in \mathbb Z\}$ in $X$. As an abelian subvariety, it is generated
by its irreducible component $Y_U^0$, containing $0$, and
by the point $\U$ of finite order in $X$, such that $U-\U\in Y_U^0,\,
N\U=\l_0\in \Lambda$.
Shifting $Z$ if needed, we may assume, without loss of generality, that
$0$ is not in the singular locus $T_U\Sigma=\Sigma\subset \Theta$.
Then $Y_U\cap \Sigma=\emptyset$, because any $T_U$-invariant set
is Zariski dense in $Y_U$. Note, that for {\it sufficiently small $t$
the affine subvariety $Y_U+Vt$ does not intersect $\Sigma$, as well}.

Consider the restriction of the theta-function onto the subspace
${\C}+Vt\subset \mathbb C^g$:
\beq\label{ttt1}
\tau (z,t)=\theta(z+Vt), \ \ z\in {\C}.
\eeq
Here and below $\C$ is a union of affine subspaces,
${\C}=\cup_{r\in \mathbb Z} (\mathbb C^d+r\U)$, where
$\mathbb C^d$ is a linear subspace that is the
irreducible component of $\pi^{-1}(Y_U^0)$, and
$\pi: \mathbb C^g\to X=\mathbb C^g/\Lambda$ is the universal cover of $X$.

The restriction of equation (\ref{cm7}) onto $Y_U$ gives the equation
\beq\label{rr}
\p_t\left(\tau(z+U,t)\,\tau(z-U,t)\right)\,\p_t\tau(z,t)=
\tau(z+U,t)\,\tau(z-U,t)\,\p^2_{tt}\tau(z,t)\ ({\rm mod}\, \tau),
\eeq
which is valid on the divisor
${\cal T}^t=\{z\in {\cal C} \mid \tau(z,t)=0\}$.
For fixed $t$, the function $u(z,t)$ has simple poles on the divisors
$\TC$ and $\TC_U=\TC-U$.

\begin{lem} Let equation (\ref{rr}) for $\tau(z,t)$ holds and let
$\l_1,\ldots,\l_d$ be a set of linear independent vectors of the
sublattice $\Lambda_U^0=\Lambda\cap \mathbb C^d\subset \mathbb C^g$. Then
for sufficiently small $t$ equation (\ref{laxd}) with the potential $u(Ux+z,t)$,
restricted to $x=n\in \mathbb Z,$ has a unique, up to a $z$-independent factor,
wave solution of  the form $\psi=k^xe^{kt}\phi(Ux+z,t,k)$
such that:

(i) the coefficients $\xi_s(z,t)$ of the formal series
\beq\label{psi2}
\phi(z,t,k)=e^{bt}\left(1+\sum_{s=1}^{\infty}\xi_s(z,t)\, k^{-s}\right)
\eeq
are meromorphic functions of the variable $z\in \mathbb C^d$
with a simple pole at the divisor $\TC$,
\beq\label{v1}
\xi_s(z,t)={\tau_s(z,t)\over \tau(z,t)}\,;
\eeq
(ii) $\phi(z,t,k)$ is quasi-periodic with respect to the lattice $\Lambda_U$
\beq\label{bloch}
\phi(z+\lambda,t,k)=\phi(z,t,k)\,B^{\,\l}(k),\ \ \lambda\in \Lambda_U;
\eeq
and is periodic with respect to the vectors $\l_0,\l_1,\ldots, \l_d$, i.e.,
\beq\label{bloch1}
B^{\,\lambda_i}(k)=1,\ \ i=0,\ldots, d.
\eeq
\end{lem}
{\it Proof.} The functions $\xi_s(z)$ are defined recursively by the equations
\beq\label{xis1}
\Delta_U\,\xi_{s+1}=\dot\xi_s+(u+b)\,\xi_s.
\eeq
Here and below $\Delta_U$ stands for the difference derivative $e^{\p_U}-1$.
The quasi-periodicity conditions (\ref{bloch}) for $\phi$ are equivalent to the equations
\beq\label{bloch2}
\xi_s(z+\lambda,t)-\xi_s(z,t)=\sum_{i=1}^s B^{\,\lambda}_i\xi_{s-i}(z,t)\,, \ \ \xi_0=1.
\eeq
A particular solution of the first equation $\Delta_U\,\xi_1=u+b$  is given by the formula
\beq\label{v5}
\xi_1^0=-\partial_V\ln \theta +l_1(z)\,b,
\eeq
where $l_1(z)$ is a linear form on ${\C}$ such that $l_1(U)=1$.
It satisfies the monodromy relations (\ref{bloch2}) with
\beq\label{bloch4}
B^{\l}_1=l_1(\l)\,b -\partial_V\ln \theta(z+\l)+\partial_V\ln \theta(z)\,,
\eeq
If $\U\neq 0$, then the space of linear forms
on ${\C}$ is $(d+1)$-dimensional. Therefore, the equation $l_1(U)=1$ and $(d+1)$
normalization conditions
$B_1^{\l_i}=1, \, i=0,1,\ldots,d,$ defines uniquely the constant $b$, the form $l_1$,
and then the constants $B_1^{\l}$ for all $\l\in \Lambda_U$.
Note, that if $\U=0$, then the space of linear forms on $\mathbb C^d$ is
$d$-dimensional, but  the normalization condition $B_1^{\l_0}$
becomes trivial.

Let us assume that the coefficient
$\xi_{s-1}$ of the series (\ref{psi2}) is known,
and that there exists a solution $\xi_s^0$ of the next equation,
which is holomorphic outside of the divisor $\TC$, and which
satisfies the quasi-periodicity conditions (\ref{bloch2}). We assume also that $\xi_s^0$
is unique up to the transformation $\xi_s=\xi_s^0+c_s(t)$, where $c_s(t)$ is
a time-dependent constant. As it was shown above,
the induction assumption holds for $s=1$.

Let us define a function $\tau_{s+1}^0(z)$ on $\TC$  with the help of
the formula
\beq\label{bl1}
\tau_{s+1}^0=-\p_t \tau_s(z,t)-b\tau_s(z,t)+{\p_t\tau(z+U,t)\over\tau(z+U,t)}\tau_s(z,t), \ \
z\in \TC.
\eeq
where $\tau_s=\theta(Ux+Vt+z)\xi_s$. Note, that the restriction of $\tau_s$
on $\TC$ does not depend on $c_s(t)$. Simple expansion of $\tau$ and
$\tau_s$ at the generic point of $\TC$ shows that the residue in $U$-line
of the right hand side of (\ref{bl1}) is equal to $R_{s+1}$, where $R_s$ is defined
by (\ref{resa}). The induction assumption (\ref{res}) of Lemma 2.2 is equivalent
to the statement: {\it if $\xi_s$ is a solution of equation (\ref{xis1}) for $s-1$, then
the function $\tau_{s+1}^0$, given by (\ref{bl1}) is equal to}
\beq\label{bl1a}
\tau_{s+1}^0=-\p_t \tau(z,t){\tau_s(z-U,t)\over\tau(z-U,t)}, \ \
z\in \TC.
\eeq
Let us show that $\tau_{s+1}^0$ is holomorphic on $\TC$.
Equations (\ref{bl1}) and (\ref{bl1a}) imply that $\tau_{s+1}^0$ is regular
on $\TC\setminus \Sigma_*$, where $\Sigma_*=\TC\cap\TC_U\cap\TC_{-U}$.
Let $z_0$ be a point of $\Sigma_*$. By the assumption,
$Y_U$ does not intersect $\Sigma$. Hence, $\TC$, for sufficiently small $t$,
does not intersect $\Sigma$, as well.
Therefore, there exists an integer $k>0$ such that $z_k=z_0-kU$ is in
$\TC$, and $\tau(z_{k+1},t)\neq 0$. Then, from equation (\ref{bl1a}) it follows that
$\tau_{s+1}^0$ is regular at the point $z=z_k$.
Using equation (\ref{bl1}) for
$z=z_k$, we get that $\p_t\tau(z_{k-1},t)\tau_s(z_k,t)=0$. The last equality
and the equation (\ref{bl1a}) for $z=z_{k-1}$ imply that
$\tau_{s+1}^0$ is regular at the point $z_{k-1}$. Regularity of
$\tau_{s+1}^0$ at $z_{k-1}$ and
equation (\ref{bl1}) for $z=z_{k-1}$ imply $\p_t\tau(z_{k-2},t)\tau_s(z_{k-1},t)=0$.
Then equation (\ref{bl1a}) for $z=z_{k-2}$ implies that
$\tau_{s+1}^0$ is regular at the point $z_{k-2}$.
By continuing these steps
we get finally that $\tau_{s+1}^0$ is regular at $z=z_0$. Therefore, $\tau_{s+1}^0$
is regular on $\TC$.

Recall, that an analytic function on an analytic divisor
in $\mathbb C^d$ has a holomorphic extension onto $\mathbb C^d$ (\cite{serr}).
The space ${\C}$ is a union of affine subspaces. Therefore,
there exists a holomorphic function $\tau^*(z,t), z\in {\C},$ such that
$\tau_{s+1}^*|_{\TC}=\tau_{s+1}^0$.
Consider the function $\chi_{s+1}^*=\tau_{s+1}^*/\tau$. It is holomorphic outside
of the divisor $\TC$. From (\ref{bloch2}) it follows that it satisfies
the relations
\beq\label{bl2}
\chi_{s+1}^*(z+\l)-\chi_{s+1}^{*}(z)=f_{s+1}^{\l}(z)+\sum_{i=1}^s B^{\,\lambda}_i
\xi_{s+1-i}(z,t)\,,
\eeq
where $f_{s+1}^{\l}(z)$ is a holomorphic function of $z\in \C$.
It satisfies the twisted homomorphism relations
\beq\label{bl3}
f_{s+1}^{\l+\mu}(z)=f_{s+1}^{\l}(z+\mu)+f_{s+1}^{\mu}(z),
\eeq
i.e., it defines an element of the first cohomology group of $\Lambda_U$ with
coefficients in the sheaf of holomorphic functions,
$f\in H^1_{gr}(\Lambda_U,H^0(\C, \O))$. The same
arguments, as that used in the proof of the part (b) of the Lemma 12 in \cite{shiota},
show that there exists a holomorphic function $h_{s+1}(z)$
such that
\beq\label{bl4}
f_{s+1}^\l(z)=h_{s+1}(z+\l)-h_{s+1}(z)+\tilde B_{s+1}^{\l},
\eeq
where $\tilde B_{s+1}^{\l}$ is a constant. Hence, the function $\chi_{s+1}=\chi^*_{s+1}+h$
has the following monodromy properties
\beq\label{bl4a}
\chi_{s+1}(z+\l)-\chi_{s+1}(z)=\tilde B_{s+1}^{\l}+\sum_{i=1}^s B^{\,\lambda}_i
\xi_{s+1-i}(z,t),
\eeq
Let us try to find a solution of (\ref{xis1}) in the form
$\xi_{s+1}^0=\chi_{s+1}+\z_{s+1}$. That gives us the equation
\beq\label{bl5}
\Delta_U\z_{s+1}=g_{s},
\eeq
where
\beq\label{bl6}
g_s=-\Delta_U\,\chi_{s+1}+\dot\xi^0_s+(u+b)\xi^0_s+\p_tc_s+(u+b)c_s.
\eeq
From (\ref{bloch2},\ref{bl4a}) it follows that $g_s$ is periodic
with respect to the lattice $\Lambda_U$, i.e., it is a function on $Y_U$. Equation
(\ref{bl2}) and the statement of Lemma 2.2 implies that it is a holomorphic function.
Therefore, $g_s$ is constant on each of the irreducible components of
$\C$:
\beq\label{apr23}
g_{\,s}(z_0+r\U)=g_{\,s}^{(r)}, \ \ z_0\in \mathbb C^d\subset \C.
\eeq
Hence, the general solution of equation (\ref{bl5}),
such that the corresponding solution $\xi_{s+1}$ of (\ref{xis1}) satisfies the
quasi-periodicity conditions, is given by the formula
\beq\label{bl7}
\z_{s+1}(z_0+r\U)=l_{s+1}(z_0)+\sum_{i=0}^{r-1} (g_{\,s}^{(i)}-a_s)+c_{s+1},
\eeq
where $l_{s+1}(z_0)$ is a linear form on $\mathbb C^d$; the constant
$a_s$ equals  $a_s=l_{s+1}(U-\U)$.
The normalization condition (\ref{bloch1}) for $B_{s+1}^{\l_i}=1, i=0,\ldots, d$
defines uniquely $l_{s+1}$ and $\p_tc_s$, i.e. the time-dependence of $c_s(t)$.
The induction step is completed and thus the lemma is proven.

Note, that a simple shift $z\to z+Z$, where $Z\notin \Sigma,$ gives
$\l$-periodic wave solutions with meromorphic coefficients along the affine
subspaces $Z+\C$. These  wave solutions are related to each other
by constant factors. Therefore choosing, in the neighborhood of any
$Z\notin \Sigma,$ a hyperplane orthogonal to the vector $U$ and
fixing initial data on this hyperplane at $t=0,$ we define the corresponding
series $\phi(z+Z,t,k)$ as a {\it local} meromorphic function of $Z$ and the
{\it global} meromorphic function of $z\in \C$.

\section{Commuting difference operators}

In this section we show that $\l$-periodic wave solutions of equation
(\ref{laxd}), constructed in the previous section, are common eigenfunctions of
rings of commuting difference operators.

\begin{lem} Let the assumptions of Theorem 1.1 hold. Then, there is a unique
pseudo-difference operator
\beq\label{LL}
\L(Z)=T+\sum_{s=0}^{\infty} w_s(Z)T^{-s}
\eeq
such that
\beq\label{kk}
\L(Ux+Vt+Z)\,\psi=k\,\psi\,,
\eeq
where $\psi=k^xe^{kt}\phi(Ux+Z,t,k)$ is a $\l$-periodic wave solution of
(\ref{laxd}).
The coefficients $w_s(Z)$ of $\L$  are meromorphic functions on the abelian variety $X$
with poles along  the divisors $T_U^{-i}\Theta= \Theta-iU,\ i\leq s$.
\end{lem}
{\it Proof.} The construction of $\L$ is standard for the theory of $2D$ Toda lattice
equations. First we define $\L$ as a pseudo-difference operator with coefficients
$w_s(Z,t)$, which are functions of $Z$ and $t$.

Let $\psi$ be a $\l$-periodic wave solution. The substitution of (\ref{psi2}) in (\ref{kk})
gives a system of equations that recursively define $w_s(Z,t)$, as difference
polynomials in the coefficients of $\psi$.
The coefficients of $\psi$ are local meromorphic functions of $Z$, but
the coefficients of $\L$ are well-defined
{\it global meromorphic functions} of on $\mathbb C^g\setminus\Sigma$, because
different $\l$-periodic wave solutions are related to each other by a factor,
which does not affect $\L$. The singular locus is
of codimension $\geq 2$. Then Hartogs' holomorphic extension theorem implies that
$w_s(Z,t)$ can be extended to a global meromorphic function on $\mathbb C^g$.

The translational invariance of $u$ implies the translational invariance of
$\L$. Indeed, for any constant $s$ the series
$\phi(Vs+Z,t-s,k)$ and $\phi(Z,t,k)$ correspond to $\l$-periodic solutions
of the same equation. Therefore, they coincide up to a $T_U$-invariant factor.
This factor does not affect $\L$. Hence, $w_s(Z,t)=w_s(Vt+Z)$.

For any $\lambda'\in \Lambda$, the $\l$-periodic wave functions corresponding to $Z$ and
$Z+\lambda'$ are also related to each other by a $T_U$-invariant factor.
Hence, $w_s$ are periodic with respect to $\Lambda$, and therefore, are
meromorphic functions on the abelian variety $X$.
The lemma is proved.

Consider now the strictly positive difference parts of the operators $\L^m$.
Let $\L^m_+$ be the difference  operator such that
$\L^m_-=\L^m-\L^m_+=F_m+F_m^1T^{-1}+O(T^{-2})$. By definition the leading
coefficient $F_m$ of $\L^m_-$ is the residue of $\L^m$:
\beq\label{res1}
F_m={\rm res}_{T}\  \L^m, \ F_m^1={\rm res}_{T}\  \L^m\,T.
\eeq
From the construction of $\L$ it follows that $[\p_t-T+u, \L^n]=0$. Hence,
\beq\label{lax}
[\p_t-T+u,\L^m_+]=-[\p_t-T+u, \L^m_-]=\left(\Delta_UF_m\right)T.
\eeq
Indeed, the left hand side of (\ref{lax}) shows that the right hand side
is a difference operator with non-vanishing coefficients only at the positive
powers of $T$. The intermediate equality shows that this operator is at most
of order $1$. Therefore, it has the form $f_m T$. The coefficient $f_m$ is
easy expressed in terms of the leading coefficient $\L^m_-$. Note, that the vanishing
of the coefficient at $T^0$ implies the equation
\beq\label{lax55}
\p_V F_m=\Delta_U\, F_m^1,
\eeq
which we will use later.

The functions $F_m(Z)$ are difference polynomials in the coefficients $w_s$ of $\L$.
Hence, $F_m(Z)$ are meromorphic functions on $X$. Next statement is crucial for
the proof of the existence of commuting differential operators associated with $u$.
\begin{lem} The abelian functions $F_m$ have at most simple poles on the divisors
$\Theta$ and $\Theta_U$.
\end{lem}
{\it Proof.} We need a few more standard constructions from $2D$ Toda theory.
If $\psi$ is as in Lemma 3.1, then  there exists a unique pseudo-difference
operator $\Phi$ such that
\beq\label{S}
\psi=\Phi k^xe^{kt},\ \ \Phi=1+\sum_{s=1}^{\infty}\f_s(Ux+Z,t)T^{-s}.
\eeq
The coefficients of $\Phi$ are universal difference polynomials in $\xi_s$.
Therefore, $\f_s(z+Z,t)$ is a global meromorphic function of $z\in \C$ and
a local meromorphic function of $Z\notin \Sigma$. Note, that $\L=\Phi\, T \Phi^{-1}$.

Consider the dual wave function defined by the left action of the operator $\Phi^{-1}$:
$\psi^+=\left(k^{-x}e^{-kt}\right)\Phi^{-1}$.
Recall that the left action of a pseudo-difference operator is the formal adjoint action
under which the left action of $T$ on a function $f$ is $(fT)=T^{-1}f$.
If $\psi$ is a formal wave solution of (\ref{laxd}),
then $\psi^+$ is a solution of the adjoint equation
\beq\label{adj}
(-\p_t-T^{-1}+u)\,\psi^+=0.
\eeq
The same arguments, as before, prove that if equations (\ref{cm5}) for poles of $v$
hold then $\xi_s^+$ have simple poles at the poles of $Tv$. Therefore, if $\psi$ as in
Lemma 2.3, then the dual wave solution is of the form
$\psi^+=k^{-x}e^{kt}\phi^+(Ux+Z,t,k)$, where
the coefficients $\xi_s^+(z+Z,t)$ of the formal series
\beq\label{psi2+}
\phi^+(z+Z,t,k)=e^{-bt}\left(1+\sum_{s=1}^{\infty}\xi^+_s(z+Z,t)\, k^{-s}\right)
\eeq
are $\l$-periodic meromorphic functions of the variable $z\in \C$ with the
simple pole at the divisor $T_U^{-1}\TC$.

The ambiguity in the definition of $\psi$ does not affect the product
\beq\label{J0}
\psi^+\psi=\left(k^{-x}e^{-kt}\Phi^{-1}\right)\left(\Phi k^xe^{kt}\right).
\eeq
Therefore, although each factor is only a local meromorphic function on
$\mathbb C^g\setminus \Sigma$, the coefficients $J_s$ of the product
\beq\label{J}
\psi^+\psi=\phi^+(Z,t,k)\,\phi(Z,t,k)=1+\sum_{s=1}^{\infty}J_s(Z,t)\,k^{-s}
\eeq
are {\it global meromorphic functions} of $Z$. Moreover, the translational invariance
of $u$ implies that they have the form $J_s(Z,t)=J_s(Z+Vt)$. The factors in
the left hand side of (\ref{J}) have the simple poles on $\Theta-Vt$ and $\Theta-U-Vt$.
Hence, $J_s(Z)$ is a meromorphic function on $X$ with the simple poles
at $\Theta$ and $T_U^{-1}\Theta=\Theta_U$.

From the definition of $\L$ it follows that
\beq\label{20}
\res_k\left(\psi^+(\L^n\psi)\right)k^{-1}dk=\res_k\left(\psi^+k^n\psi\right)k^{-1}dk
=J_{n}.
\eeq
On the other hand, using the identity
\beq\label{dic}
\res_k \left(k^{-x}\D_1\right)\left(\D_2k^x\right)k^{-1}dk=
\res_{T}\left(\D_2\D_1\right),
\eeq
we get
\beq\label{201}
\res_k(\psi^+\L^n\psi)k^{-1}dk=
\res_k\left(k^{-x}\Phi^{-1}\right)\left(\L^n\Phi k^x\right)k^{-1}dk=
\res_{T}\L^n=F_n.
\eeq
Therefore, $F_n=J_{n}$ and the lemma is proved.

Let $\bf {\hat F}$ be a linear space generated by $\{F_m, \, m= 1,\ldots\}$.
It is a subspace of the $2^g$-dimensional space of the abelian functions that have
at most simple poles at
$\Theta$ and $\Theta_U$. Therefore, for all but $\hat g=\dim\ {\bf \hat F}$ positive integers $n$,
there exist constants $c_{i,n}$ such that
\beq\label{f1}
F_n(Z)+\sum_{i=1}^{n-1} c_{i,n}F_i(Z)=0.
\eeq
Let $I$ denote the subset of integers $n$ for which there are no such constants. We call
this subset the gap sequence.
\begin{lem} Let $\L$ be the pseudo-difference operator corresponding to
a $\l$-periodic wave function $\psi$ constructed above. Then, for the difference
operators
\beq\label{a2}
L_n=\L^n_++\sum_{i=1}^{n-1} c_{i,n}\L^{n-i}_+=0, \ n\notin I,
\eeq
the equations
\beq\label{lp}
L_n\,\psi=a_n(k)\,\psi, \ \ \ a_n(k)=k^n+\sum_{s=1}^{\infty}a_{s,n}k^{n-s}
\eeq
where $a_{s,n}$ are constants, hold.
\end{lem}
{\it Proof.} First note that from (\ref{lax}) it follows that
\beq\label{lax3}
[\p_t-T+u,L_n]=0.
\eeq
Hence, if $\psi$ is a $\l$-periodic wave solution of (\ref{laxd})
corresponding to $Z\notin \Sigma$, then $L_n\psi$ is also a $\l$-periodic
solution of the same equation. That implies the equation
$L_n\psi=a_n(Z,k)\psi$, where $a$ is $T_U$-invariant. The ambiguity in the definition of
$\psi$ does not affect $a_n$. Therefore, the coefficients of $a_n$ are well-defined
{\it global} meromorphic functions on $\mathbb C^g\setminus \Sigma$. The $\p_U$-
invariance of $a_n$ implies that $a_n$, as a function of $Z$, is holomorphic outside
of the locus. Hence it has an extension to a holomorphic function on $\mathbb C^g$.
It is periodic with respect to the lattice
$\Lambda$. Hence $a_n$ is $Z$-independent. Note that $a_{s,n}=c_{s,n},\ s\leq n$.
The lemma is proved.

The operator $L_m$ restricted to the points $x=n$
can be regarded as a ${Z}$-parametric family   of
ordinary difference operators $L_m^{Z}$, whose coefficients have the form
\beq\label{lu}
L_m^Z=T^m+\sum_{i=1}^{m-1} u_{i,m}(Un+Z)\, T^{m-i},\ \ m\notin I.
\eeq
\begin{cor} The operators $L_m^{Z}$
commute with each other,
\beq\label{com1}
[L_n^{Z},L_m^{Z}]=0, \ \ Z\notin \Sigma.
\eeq
\end{cor}
From (\ref{lp}) it follows that $[L_n^{Z},L_m^{Z}]\psi=0$.
The commutator is an ordinary difference operator.
Hence, the last equation implies (\ref{com1}).

\section{The spectral curve.}

A theory of commuting difference operators containing a pair of operators of co-prime orders
was developed in  (\cite{mum,kr-dif}). It is analogous to the theory of rank 1 commuting
differential operators (\cite{ch1,ch2,kr1,kr2,mum}). (Relatively recently
this theory was generalized to the case of commuting difference operators of arbitrary rank
in \cite{n-kr}.)

\begin{lem} Let $\A^Z,\ Z\notin \Sigma,$ be a commutative ring of ordinary difference
operators spanned by the operators $L_n^Z$. Then there is an irreducible algebraic
curve $\G$ of arithmetic genus $\hat g=\dim\ {\bf \hat F}$, such that for a generic $Z$
the ring $\A^Z$ is isomorphic to the ring $A(\G,P_+,P_-)$ of the meromorphic functions
on $\G$ with the only pole at a smooth point $P_+$, vanishing at another smooth point
$P_-$. The correspondence $Z\to \A^Z$ defines a holomorphic map of
$X\setminus \Sigma$ to the space of torsion-free rank 1 sheaves $\F$ on $\G$
\beq\label{is}
j: X\backslash\Sigma\longmapsto \overline{\rm Pic}(\G).
\eeq
\end{lem}
{\it Proof.} As shown in (\cite{mum,kr-dif})
there is a natural correspondence
\beq\label{corr}
\A\longleftrightarrow \{\G,P_{\pm},  \F\}
\eeq
between commutative rings $\A$ of ordinary linear
difference operators containing a pair of monic operators of co-prime orders, and
sets of algebraic-geometrical data $\{\G,P_{\pm}, [k^{-1}]_1, \F\}$, where $\G$ is an
algebraic curve with a fixed first jet $[k^{-1}]_1$ of a local coordinate $k^{-1}$ in the neighborhood of a smooth
point $P_+\in\G$ and $\F$ is a torsion-free rank 1 sheaf on $\G$ such that
\beq\label{sheaf}
h^0(\G,\F(nP_+-nP_-))=h^1(\G,\F(nP_+-nP_-)=0.
\eeq
The correspondence becomes one-to-one if the rings $\A$ are considered modulo conjugation
$\A'=g(x)\A g^{-1}(x)$.

The construction of the correspondence (\ref{corr})
depends on a choice of initial point $x_0=0$. The spectral curve and the sheaf $\F$
are defined by the evaluations of the coefficients of generators of $\A$
at a finite number of points of the form $x_0+n$.
In fact, the spectral curve is independent on the choice of $x_0$, but the sheaf
does depend on it, i.e. $\F=\F_{x_0}$.

Using the shift of the initial point it is easy to show that the correspondence
(\ref{corr}) extends to the commutative rings of operators whose coefficients are
{\it meromorphic} functions of $x$. The rings of operators having poles at $x=0$
correspond to sheaves for which the condition (\ref{sheaf}) for $n=0$
is violated.

The algebraic curve $\G$ is called the spectral curve of $\A$.
The ring $\A$ is isomorphic to the ring $A(\G,P_+,P_-)$ of meromorphic functions
on $\G$ with the only pole at the puncture $P_+$ and which vanish at $P_-$.
The isomorphism is defined by
the equation
\beq\label{z2}
L_a\psi_0=a\psi_0, \ \ L_a\in \A, \ a\in A(\G,P_+,P_-).
\eeq
Here $\psi_0$ is a common eigenfunction of the commuting operators. At $x=0$ it is
a section of the sheaf $\F\otimes\O(P_+)$.

Let $\G^Z$ be the spectral curve corresponding to $\A^Z$.
It is well-defined for all $Z\notin \Sigma$.
The eigenvalues $a_n(k)$ of the operators $L_n^Z$ defined in (\ref{lp})
coincide with the Laurent expansions at $P_+$ of the meromorphic
functions $a_n\in A(\G^Z,P_{\pm})$. They are $Z$-independent. Hence, the spectral curve
is $Z$-independent, as well, $\G=\G^Z$. The first statement of the lemma is thus
proven.

The construction of the correspondence
(\ref{corr}) implies that if the coefficients of operators $\A$
holomorphically depend on parameters then the algebraic-geometrical spectral data are
also holomorphic functions of the parameters. Hence $j$ is holomorphic away of $\Theta$.
Then using the shift of the initial point and the fact, that $\F_{x_0}$ holomorphically
depends on $x_0$, we get that $j$ holomorphically extends on $\Theta\setminus \Sigma$,
as well. The lemma is proved.

{\it Remark.}
Recall, that a commutative ring $\A$ of linear ordinary difference operators
is called maximal if it is not contained in any bigger commutative ring.
As in the differential case (see details in \cite{kr-schot}), it is
easy to show that for the generic $Z$ {\it the ring $\A^Z$ is maximal}.

Our next goal is to prove finally the global existence of the wave function.
\begin{lem} Let the assumptions of the Theorem 1.1 hold. Then there exists a common
eigenfunction of the operators $L_n^Z$ of the form
$\psi=e^{kx}\phi(Ux+Z,k)$ such that
the coefficients of the formal series
\beq\label{psi6}
\phi(Z,k)=1+\sum_{s=1}^{\infty}\xi_s(Z)\, k^{-s}
\eeq
are global meromorphic functions with a simple pole at $\Theta$.
\end{lem}
{\it Proof.} It is instructive to consider first the case when the spectral curve $\G$
of the rings $\A^Z$ is smooth. Then, as shown in (\cite{kr-dif}),
the corresponding common
eigenfunction of the commuting differential operators (the Baker-Akhiezer function),
normalized by the
condition $\psi_0|_{x=0}=1$, is of the form
\beq\label{ba}
\hat \psi_{0}={\hat\theta (\hat A(P)+\hat Ux+\hat Z)\,\hat\theta (\hat Z)\over
\hat\theta(\hat Ux+\hat Z)\,\hat\theta(\hat A(P)+\hat Z)}\,
e^{x\,\Omega(P)}.
\eeq
Here $\hat \theta (\hat Z)$ is the Riemann theta-function constructed with the help of the
matrix of $b$-periods of normalized holomorphic differentials on $\G$;
$\hat A: \G\to J(\G)$is the Abel map; $\Omega$ is the abelian integral corresponding to the third kind
meromorphic differential $d\Omega$ with the residues $\pm 1$ at the punctures
$P_{\pm}$ and $2\pi i \hat U$ is the vector of its $b$-periods.

\noindent{\it Remark.}
Let us emphasize, that the formula (\ref{ba}) is not the result of a solution of some
difference equations. It is a direct corollary of analytic properties of
the Baker-Akhiezer function $\hat \psi_0(x,P)$ on the spectral curve:

$(i)$ {\it $\hat\psi_0$ is meromorphic function on the universal
cover $\tilde \G$ of $\{\G\setminus P_{\pm}\}$ with the monodromy around $P_{\pm}$ equals
$e^{\pm 2\pi i x}$; the pole divisor od $\hat \psi_0$
is of degree $\tilde g$ and is $x$-independent. It is non-special,
if the operators are regular at the normalization point $x=0$};

$(ii)$ {\it in the neighborhood of $P_0$ the function $\hat\psi_0$ has the form (\ref{ps})
(with $t=0$)}.

\noindent
From the  Riemann-Rokh theorem it follows that, if $\hat\psi_0$ exists, then it is unique.
It is easy to check that the function $\hat\psi_0$ given by (\ref{ba})
has all the desired properties.

The last factors in the numerator and the denominator of (\ref{ba}) are $x$-independent.
Therefore, the function
\beq\label{ba1}
\hat \psi_{BA}={\hat\theta (\hat A(P)+\hat Ux+\hat Z)\over
\hat\theta(\hat Ux+\hat Z)}\,
e^{x\,\Omega(P)}
\eeq
is also a common eigenfunction of the commuting operators.

In the neighborhood of $P_+$ the function $\hat \psi_{BA}$ has the form
\beq\label{sh10}
\hat \psi_{BA}=k^x
\left(1+\sum_{s=1}^{\infty}{\tau_s (\hat Z+\hat Ux)\over \hat \theta(\hat U x+\hat Z)}\,k^{-s}
\right), \ \ k=e^{\Omega},
\eeq
where $\tau_s(\hat Z)$ are global holomorphic functions.

According to Lemma 4.1, we have a holomorphic map $\hat Z=j(Z)$
of $X\setminus\Sigma$ into $J(\G)$.
Consider the formal series $\psi=j^*\hat \psi_{BA}$.
It is globally well-defined out of $\Sigma$.
If $Z\notin \Theta$, then $j(Z)\notin \hat \Theta$
(which is the divisor on which the condition (\ref{sheaf}) is violated).
Hence, the coefficients of $\psi$ are regular out of $\Theta$. The singular locus
is at least of codimension 2. Hence, using once again Hartogs' arguments we can
extend $\psi$ on $X$.

If the spectral curve is singular, we can proceed along the same lines
using a proper generalization of (\ref{ba1}). Note, that in (\cite{kr-schot})
we used the generalization of (\ref{ba1}) given by the theory of Sato $\tau$-function
(\cite{wilson}). In fact the general theory of the tau-function is not needed for
our purposes. It is enough to consider only algebro-geometric points
of the Sato Grassmanian.

Let $p:\hat \G\longmapsto \G$ be the normalization map, i.e. a regular map
of a smooth genus $\tilde g$ algebraic curve $\hat \G$ to the spectral curve
$\G$ which is one-to-one outside the preimages $\hat q_k$ of singular points $\G$.
The normalized common eigenfunction $\hat \psi_0$ of commuting operators can be regarded
as a multi-valued meromorphic function on ${\hat \G\setminus P_{\pm}}$ with the monodromy
$e^{\pm 2\pi i x}$ around the punctures $P_{\pm}$. The divisor $D=\sum_s \g_s$ of
the poles of $\hat\psi_0$ is of degree $\tilde g+d\leq \hat g$, where $\hat g$ is the
arithmetic genus of $\G$. The expansions of $\hat\psi$ at the points $\hat q_k$ are in
some linear subspace of co-dimension $d$ in the space $\bigoplus_{k}\O_{\hat q_k}$. If
we fix local coordinates $z_k$ in the neighborhoods of $\hat q_k$, then the latter
condition can be written as a system of $n$ linear constraints
\beq\label{constr}
\sum_{k,j}c_{k,s}^i\p_{z_k}^s\hat\psi_0|_{q_k}=0, \ \ i=1,\ldots,d.
\eeq
The constants $\hat c_{k,j}^i$ are defined up to the
transformations $\hat c_{k,s}^i\to \sum_i g_i^{i'}\hat c_{k,s}^i$,
where $g_i^{i'}$ is a non-degenerate matrix.

The analytical properties of the function $\hat \psi_0$ imply that it can be represented
in the form
\beq\label{bl10}
\hat \psi_0=\sum_{i=0}^d r_i(x,D) {\hat\theta (\hat A(P)+\hat Ux+\hat Z_i)\over
\hat\theta(\hat A(P)+\hat Z_i)}\,
e^{x\,\Omega(P)}.
\eeq
Here $Z_i=Z_*-\hat A(\g_{\tilde g+i})$, where
$Z_*=R-\sum_{s=1}^{\tilde g-1}\hat A(\g_s)$ and $R$ is the vector of the Riemann
constants.

The coefficients $r_i$ in (\ref{bl10}) are defined by the linear equations (\ref{constr})
and the normalization of the leading term in the expansion (\ref{ps})
of $\hat \psi_0$ at $P_+$. Keeping track of the $x$-dependent terms one
can write the eigenfunction of the commuting operators in the
form
\beq\label{bl11}
\hat\psi_{BA}=\sum_{i=0}^n R_i \, \hat\theta (\hat A(P)+\hat Ux+\hat Z_i),
\eeq
where the coefficients $R_i$ depend on $x$, but are $P$-independent. The equations
(\ref{constr}) imply
\beq\label{bl12}
\sum_{j=0}^n M_{i,j}R_j=0,
\eeq
where the entries of $(d\times (d+1))$-dimensional matrix $M_{ij}$ are equal to
\beq\label{bl13}
M_{ij}=\sum_{k,s}C_{k,s}^{\,i}\p_{z_k}^s \hat\theta (\hat A(P(z_k))+\hat Ux+\hat Z_j)|_{z_k=0}.
\eeq
The coefficients $C_{k,s}^j$ in (\ref{bl13}) can be expressed in terms of $D$
and the coefficients $c_{k,s}^i$ in (\ref{constr}). They and the divisor $D$ can be
regarded as parameters defining the sheaf $\F$ in (\ref{corr}).

Let us define a function $\tau(x,P;\F)$ as the determinant of
$(d+1)\times (d+1)$-dimensional matrix
\beq\label{bl14}
\tau(x,P;\F)=\det\left(\begin{array}{ccc}
\hat\theta (\hat A(P)+\hat Ux+\hat Z_0)   &\cdots        &
\hat\theta (\hat A(P)+\hat Ux+\hat Z_d)\\
M_{1,0}      &\cdots        & M_{1,d}\\
\cdot        &  \cdot    & \cdot    \\
M_{n,0}  &  \cdots         &M_{d,d}
\end{array}\right).
\eeq
Then, the common eigenfunction of the commuting operators
can be represented in the form
\beq\label{baw}
\hat \psi_{BA}={\tau(x,P\,; \F) \over
\tau(x,P_+; \F)}e^{x\Omega(P)}
\eeq
The rest of the arguments proving the lemma are
the same, as in the smooth case.

\begin{lem} There exist $g$-dimensional vectors $V_{m}=\{V_{m,k}\}$
and constants $v_{m}$ such that the abelian functions $F_m=\res_T \L^m$ are equal to
\beq\label{nnov7}
F_{m}(Z)=v_{m}+ \Delta_U\left(\p_{\,V_m}\ln \theta(Z)\right),
\eeq
where $\p_{\,V_m}=\sum_{k=1}^g V_{m,k}\p_{z_k}$.
\end{lem}
{\it Proof.} The proof is identical to that of Lemma 3.6 in \cite{kr-schot}.
Recall that the functions $F_n$ are abelian
functions with simple poles at the divisors $\Theta$ and $\Theta_U$.
In order to prove the statement of the lemma
it is enough to show that $F_n=\Delta_U Q_n$, where $Q_n$ is a
meromorphic function with a pole along $\Theta$. Indeed, if $Q_n$ exists,
then, for any vector $\l$ in the period lattice, we have $Q_n(Z+\l)=Q_n(Z)+c_{n,\l}$.
There is no abelian function with a simple pole on $\Theta$. Hence, there exists
a constant $q_n$ and two $g$-dimensional vectors $l_n,V_n$, such that
$Q_n=q_n+(l_n,Z)+(V_n,h(Z))$, where $h(Z)$ is a vector with the coordinates
$\p_{z_i}\ln\theta$. Therefore, $F_n=(l_n,U)+(V_n,\Delta_U)h$.

Let $\psi(x,Z,k)$ be the formal Baker-Akhiezer function defined in the previous lemma.
Then the coefficients $\varphi_s(Z)$ of the corresponding wave operator $\Phi$
are global meromorphic functions with poles on $\Theta$.

The left and  right action of pseudo-difference operators are formally adjoint,
i.e., for any two operators the equality
$\left(k^{-x}\D_1\right)\left(\D_2k^{x}\right)=
k^{-x}\left(\D_1\D_2k^x\right)+(T-1)\left(k^{-x}\left(\D_3k^x\right)\right)$
holds. Here $\D_3$ is a pseudo-difference operator whose coefficients are difference
polynomials in the coefficients of $\D_1$ and $\D_2$. Therefore, from (\ref{J0}-\ref{201})
it follows that
\beq\label{z8}
\psi^+\psi=1+\sum_{s=2}^{\infty}F_{s-1}k^{-s}=
1+\Delta\left(\sum_{s=2}^{\infty}Q_sk^{-s}\right).
\eeq
The coefficients of the series $Q$ are difference polynomials in the
coefficients $\varphi_s$ of the wave operator.
Therefore, they are global meromorphic
functions of $Z$ with poles on $\Theta$. Lemma is proved.

In order to complete the proof of our main result we need one more standard fact
of the $2D$ Toda lattice theory: flows of the $2D$ Toda lattice
hierarchy define deformations of the commutative
rings $\A$ of ordinary linear difference operators.
The spectral curve is invariant under these flows. There are two sets
of $2D$ Toda hierarchy flows. Each of them is isomorphic to
the KP hierarchy.
For a given spectral curve $\G$ the orbits of
the KP hierarchy are isomorphic to the generalized Jacobian $J(\G)={\rm Pic}^0 (\G)$,
which is the equivalence classes of zero degree divisors on the spectral curve.
(see \cite{shiota,zab,kr-dif,wilson}).

The part of $2D$ hierarchy we are going to use
is a system of  commuting differential
equation for a pseudo-difference operator $\L$
\beq\label{z4}
\p_{t_n}\L=[\L^n_++F_n,\L]=-[\L_-^n-F_n,\L]\,.
\eeq
The coefficient $w_0$ of $\L$ in (\ref{LL}) equals
\beq\label{kkk}
w_0=-\Delta_U\xi_1=-u.
\eeq
Therefore, (\ref{z4}) and equations (\ref{lax55}, \ref{nnov7}) imply
\beq\label{z5}
\p_{t_n}u=-\Delta_U F_n^1=-\p_VF_n=-\p_V [\Delta_U\p_{V_n}\ln \theta(Z)\,,
\eeq
where $V_n$ is the vector defined in (\ref{nnov7}).

Equations (\ref{z5}) identify the tangent vector $\p_{t_n}$ to the orbit of
the KP hierarchy with the tangent vector $\p_{V_n}$ to the abelian variety $X$.
Hence, for a generic $Z\notin \Sigma$, the orbit of the KP flows
of the ring $\A^Z$ is in $X$, i.e. it defines an holomorphic imbedding:
\beq\label{imb}
i_Z:J(\G)\longmapsto X.
\eeq
From (\ref{imb}) it follows that $J(\G)$ is {\it compact}.

The generalized Jacobian of an algebraic curve is compact
if and only if the curve is {\it smooth} (\cite{mdl}). On a smooth algebraic curve
a torsion-free rank 1 sheaf is a line bundle, i.e. $\overline {\rm Pic} (\G)=J(\G)$.
Then (\ref{is}) and the dimension arguments imply that $i_Z$ is an isomorphism
and the map  $j$ is  inverse to $i_Z$. Theorem 1.1 is proved.

\section{Fully discrete case}

In this section we present the proof of Theorem 1.2. As above, we begin with
the proof of the implication $(A) \longmapsto (C)$. We would like to mention
that equation (\ref{cm7d}) can be derived as a necessary condition for the existence of
a solution of (\ref{laxdd}), which is meromorphic in any of the variables $m,n$ or in
their linear combinations. For further use, let us introduce the variables
\beq\label{var}
x=m-n,\ \ \nu=m+n-1\,.
\eeq
In these variables equation (\ref{laxdd}) takes the form
\beq\label{laxm}
\psi(x-1,\nu)=\psi(x+1,\nu)+u(x,\nu)\psi(x,\nu-1)\,.
\eeq
Let $\tau(x,\nu)$ be a holomorphic function of the variable $x$ in some {\it
translational invariant} domain $\D=T\D\in \mathbb C$, where $T:x\to x+1$.
Suppose that for each $\nu$ the function $\tau$ in $\D$ has a simple root $\eta(\nu)$
such that
\beq\label{xid}
\tau(\eta(\nu)+1,\nu-1)\,\tau(\eta(\nu)-1,\nu-1)\neq 0.
\eeq
\begin{lem}
If equation (\ref{laxm}) with the potential
\beq\label{f1d}
u={\tau(x,\nu+1)\,\tau(x,\nu-1)\over
\tau(x-1,\nu)\,\tau(x+1,\nu)}
\eeq
has a meromorphic in $\D$ solution $\psi(x,\nu)$ such that it has a
simple pole at $\eta(\nu)$, and regular at $\eta(\nu+1)-1, \eta(\nu+1)+1$, then the equation
\beq\label{f2d}
{\tau(\eta+1,\nu+1)\,\tau(\eta-2,\nu)\,\tau(\eta+1,\nu-1)\over
\tau(\eta-1,\nu+1)\,\tau(\eta+2,\nu)\,\tau(\eta-1,\nu-1)}=-1\,,\ \ \eta=\eta(\nu)
\eeq
holds.
\end{lem}
{\it Proof.} The substitution in (\ref{laxm}) of the Laurent expansion (\ref{psie})
for $\psi$ (with coefficients depending on $\nu$), and
the expansion
\beq\label{ued}
\tau(x,\nu)=v_0(\nu)\,(x-\eta(\nu))+O((x-\eta(\nu))^2)\,,
\eeq
gives the following system of equations.

From the vanishing of the residues at $\eta+1$ and $\eta-1$
of the left hand side of (\ref{laxm}) we get
\begin{eqnarray}
\a(\nu)&=&{\tau(\eta+1,\nu+1)\,\tau(\eta+1,\nu-1)\over
\tau(\eta+2,\nu)\,v_0(\nu)}\,\psi(\eta+1,\nu-1), \label{eq1d}\\
-\a(\nu)&=&{\tau(\eta-1,\nu+1)\,\tau(\eta-1,\nu-1)\over
\tau(\eta-2,\nu)\,v_0(\nu)}\, \psi(\eta-1,\nu-1), \label{eq2d}
\end{eqnarray}
The evaluation of (\ref{laxm}) at $x=\eta(\nu+1)$ gives the equation
\beq\label{eq3d}
\psi(\eta(\nu+1)-1,\nu)=\psi(\eta(\nu+1)+1,\nu).
\eeq
Equations (\ref{eq1d},\ref{eq2d},\ref{eq3d}) directly imply (\ref{f2d}). The Lemma is
proved.

Equation (\ref{f2d}) for $\tau$ of the form $\theta(mU+nV+Z)$ coincides with
(\ref{cm7d}), i.e., the implication $(A) \longmapsto (C)$ is proved.

The next step is to show that equations (\ref{f2d}) are sufficient for local existence
of wave solutions with coefficients having poles only at zeros of $\tau$. The wave
solutions of (\ref{laxm}) are solutions of the form
\beq\label{psdd}
\psi(x,\nu,k)=k^{\nu}
\left(1+\sum_{s=1}^{\infty}\xi_s(x,\nu)\,k^{-s}\right)\,.
\eeq
\begin{lem} Suppose that $\tau(x,\nu)$ is holomorphic in a domain $\D$ of the form
(\ref{dd}) where it has simple zeros, for which condition (\ref{xid}) and
equation (\ref{f2d}) hold.
Then there exist meromorphic wave solutions of equation (\ref{laxm}) that have
simple poles at zeros of $\tau$ and are holomorphic everywhere else.
\end{lem}
{\it Proof.} Substitution of (\ref{psdd}) into (\ref{laxm}) gives a recurrent system of
equations
\beq\label{xisdd}
\xi_{s+1}(x-1,\nu)-\xi_{s+1}(x+1,\nu)=u(x,\nu)\,\xi_s(x,\nu-1).
\eeq
Under the assumption that $\D$ is a disconnected union of small disks, $\xi_{s+1}$
can be defined as an arbitrary meromorphic function in $D_0$ and then
extended on $\D$ with the help of (\ref{xisdd}).  Our goal is to prove
by induction that (\ref{xisdd}) has meromorphic solutions with simple poles only at the zeros
of $\tau$.

Suppose that $\xi_s(x,\nu)$ has a simple pole at $x=\eta(\nu)=\eta$
\beq\label{5dd}
\xi_s={r_s\over x-\eta}+r_{s0}+\cdots\,.
\eeq
The condition that $\xi_{s+1}(x,\nu)$ has no pole at $\eta+1$ is equivalent to the equation
\beq\label{resdd}
r_{s+1}(\nu)={\tau(\eta+1,\nu+1)\,\tau(\eta+1,\nu-1)\over
v_0\,\tau(\eta+2,\nu)}\ \xi_s(\eta+1,\nu-1).
\eeq
The condition that $\xi_{s+1}(x,\nu)$ has no pole at $\eta-1$
is equivalent to the equation
\beq\label{resdd1}
-r_{s+1}(\nu)={\tau(\eta-1,\nu+1)\,\tau(\eta-1,\nu-1)\over
v_0\,\tau(\eta-2,\nu)}\ \xi_s(\eta-1,\nu-1).
\eeq
Using (\ref{xisdd}) for $s-1$, and the equation $u(\eta,\nu-1)=0$, we get
\beq\label{fd4d}
\xi_s(\eta-1,\nu-1)=\xi_s(\eta+1,\nu-1).
\eeq
Then,  equation (\ref{f2d}) imply that two
different expressions for $r_{s+1}(\nu)$ obtained from (\ref{resdd}) and (\ref{resdd1})
do in fact coincide. The lemma is proved.

\bigskip
\noindent
{\bf Normalization problem.}
As before, our goal is to show that wave solutions can be defined uniquely
up to a constant factor with the help of certain quasi-periodicity conditions.
That requires the global existence of the wave functions along certain
affine subspaces. In what follows we  use affine subspaces in the direction of the vector
$2W=U-V$. The {\it singular locus} $\Sigma$ which controls
the obstruction for the global existence of the wave solution on $X$ is
the maximal $T_{U-V}$-invariant subset of the theta-divisor, and which is not
$T_U$ or $T_V$-invaiant
\beq\label{locus}
\Sigma=\left\{Z\ \Big|\
{\theta (k(U-V)+Z)\over \theta(Z+U)}=0,
\ \ {k(U-V)+Z)\over \theta(Z+V)}=0\ \ \ \forall \, k\in \mathbb Z\,\right\}\,.
\eeq
Surprisingly it turns out that in the fully discrete case the proof
of the statement that the singular locus is in fact empty can be obtained at
much earlier stage than in the continuous or semi-continuous cases.

Let $Y=\langle (U-V)\,k\rangle$ be the Zariski closure
of the group $\{(U-V)\, k \mid k\in \mathbb Z\}$ in $X$. It is generated
by its irreducible component $Y^0$, containing $0$, and
by the point $W_0$ of finite order in $X$, such that $2W-W_0\in Y^0,\,
NW_0=\l_0\in \Lambda$. Shifting $Z$ if needed, we may assume, without loss of
generality, that $0$ is not in the singular locus. Then $Y\cap \Sigma=\emptyset$.

Let $\tau(z,\nu)$ be a function defined by the formula
\beq\label{td}
\tau (z,\nu)=\theta\left(z+\frac{\nu}{2}\,(U+V)\right), \ \ z\in {\C}.
\eeq
Here and below $\C$ is a union of affine subspaces, that are preimages
of irreducible components of $Y$ under the projection
$\pi: \mathbb C^g\to X=\mathbb C^g/\Lambda$.

The restriction of equation (\ref{cm7d}) onto $Y$ gives the equation
\beq\label{f5d}
{\tau(z+W,\nu+1)\,\tau(z-2W,\nu)\,\tau(z+W,\nu-1)\over
\tau(z-W,\nu+1)\,\tau(z+2W,\nu)\,\tau(z-W,\nu-1)}=-1\,,
\eeq
which is valid on the divisor ${\cal T}^{\nu}=\{\,z\in{\cal C}\,\mid\, \tau(z,\nu)=0\}$.

The function
\beq\label{f7d}
u={\tau(z,\nu+1)\,\tau(z,\nu-1)\over
\tau(z-W,\nu)\,\tau(z+W,\nu)}
\eeq
is periodic with respect to the lattice $\Lambda_W=\Lambda\cap {\C}$.
The latter is generated by the sublattice $\Lambda_W^0=\Lambda\cap \mathbb C^d$,
where $\mathbb C^d$ is a linear subspace in $\mathbb C^g$, that is preimage
of $Y_W^0$, and the vector $\l_0=NW_0\in \Lambda$.
For fixed $\nu$, the function $u(z,t)$ has simple poles on the divisors
$\TD\pm W$.

\begin{lem} Let $\tau(z,\nu)$ be a sequence of non-trivial holomorphic functions
on $\C$ such that $u(z,\nu)$ given by (\ref{f7d}) is periodic with respect to $\Lambda_W$.
Suppose that equation (\ref{f5d}) holds.
Then  there exist wave solutions $\psi(z,\nu,k)=k^\nu\phi(z,\nu,k)$
of the equation
\beq\label{laxm1}
\psi(z-W,\nu,k)=\psi(z+W,\nu,k)+u(z,\nu)\psi(z,\nu-1,k)\,,
\eeq
such that:

(i) the coefficients $\xi_s(z,\nu)$ of the formal series
\beq\label{psi2d}
\phi(z,\nu,k)=1+\sum_{s=1}^{\infty}\xi_s(z,\nu)\, k^{-s},
\eeq
are meromorphic functions of the variable $z\in \mathbb \C$ with simple poles at the
divisor $\TD$, i.e.
\beq\label{new}
\xi_s(z,\nu)={\tau_s(z,\nu)\over \tau(z,\nu)},
\eeq
where $\tau_s(z,\nu)$ is now a holomorphic function;

(ii) $\xi_s(z,\nu)$ satisfy the following monodromy properties
\beq\label{new1}
\xi_s(z+\lambda,\nu)-\xi_{s}(z,\nu)=\sum_{i=1}^s B^{\,\lambda}_{i,\,\nu-s+i}\,\xi_{s-i}(z,\nu)\,,\ \
\ \lambda\in \Lambda_W,
\eeq
where $B^{\,\lambda}_{i,\,\nu}$ are $z$-independent.
\end{lem}
{\it Proof.} The functions $\xi_s(z,\nu)$ are defined recursively by the equations
\beq\label{laxm2}
\xi_{s+1}(z-W,\nu)-\xi_{s+1}(z+W,\nu)=u(z,\nu)\,\xi_s(z,\nu-1).
\eeq
We will now prove lemma by induction in $s$. Let us assume inductively that that for
$r\leq s$ the functions $\xi_r$ are known and satisfy (\ref{new1}).
Then, we define the residue of $\xi_{s+1}$ on $\TD$ by formulae
\begin{eqnarray}
\tau_{s+1}^0(z,\nu)&=&{\tau(z+W,\nu+1)\,\tau_s(z+W,\nu-1)\over \tau(z+2W,\nu)}\,,\ \
z\in \TD,\label{may2}\\
-\tau_{s+1}^0(z,\nu)&=&{\tau(z-W,\nu+1)\,\tau_s(z-W,\nu-1)\over \tau(z-2W,\nu)}\,,\ \
z\in \TD\,, \label{may3}
\end{eqnarray}
which, as follows from (\ref{resdd}) and (\ref{resdd1}), coincide.
The expression (\ref{may2}) is certainly holomorphic when
$\tau(z+2W)$ is non-zero, i.e. is holomorphic
outside of $\TD\cap(\TD-2W)$. Similarly from
(\ref{may3}) we see that $\tau_{s+1}^0(z,\nu)$ is
holomorphic away from  $\TD\cap(\TD+2W)$.

We claim that $\tau_{s+1}^0(z,\nu)$ is holomorphic
everywhere on $\TD$. Indeed, by definition of $Y$,
the closure of the abelian subgroup generated
by $2W$ is everywhere dense. Thus for any $z\in\TD$
there must exist some $N\in\mathbb N$ such that
$z-2(N+1)W\not\in\TD$; let $N$ moreover be the
minimal such $N$. From (\ref{may3}) it then follows that
$\tau_{s+1}^0(z,\nu)$ can be extended holomorphically to the point
$z-2NW$. Thus expression
(\ref{may2}) must also be holomorphic at $z-2NW$; since its
denominator there vanishes, it means that the numerator must also
vanish.
But this expression is equal to the numerator of (\ref{may3}) at
$z-2(N-1)W$; thus $\tau_{s+1}^0$ defined from
(\ref{may3}) is also holomorphic at $z-2(N-1)W$ (the
numerator vanishes, and the vanishing order of the denominator is
one, since we are talking exactly about points on its vanishing
divisor).
Note that we did not quite need the fact $z-2(N+1)W\not\in\TD$
itself, but the consequences of the minimality of $N$, i.e.,
$z-2kW\in\TD$, $0\leq k\leq N$,
and the holomorphicity of $\tau_{s+1}^0(z,\nu)$ at $z-2NW$."
Therefore, in the same
way, by replacing $N$ by $N-1$, we can then deduce holomorphicity
$\tau_{s+1}^0(z,\nu)$ at $z-2(N-2)W$ and, repeating the process $N$ times, at
$z$.

Recall once again that that an analytic function on an analytic divisor
in $\mathbb C^d$ has a holomorphic extension to all of $\mathbb C^d$ (\cite{serr}).
Therefore, there exists a holomorphic function $\widetilde \tau_{s+1}(z,\nu)$
extending the $\tau_{s+1}^0(z,\nu)$. Consider then the function
$\chi_{s+1}(z,\nu)=\wt\tau_{s+1}(z,\nu)/\tau (z,\nu)$,
holomorphic outside of $\TD$ .

From (\ref{new1}) and (\ref{may2}) it follows that the function
\beq\label{new6}
f_{s+1}^{\l}(z,\nu)=\chi_{s+1}(z+\l,\nu)-
\chi_{s+1}(z,\nu)-\sum_{i=1}^{s} B^{\,\lambda}_{i,\,\nu-1-s+i}\,
\xi_{s+1-i}(z,\nu)
\eeq
vanishes at the divisor $\TD$. Hence, it  is a holomorphic function.
It satisfies the twisted homomorphism relations
\beq\label{new3}
f_{s+1}^{\l+\mu}(z,\nu)=f_{s+1}^{\l}(z+\mu,\nu)+f_{s+1}^{\mu}(z,\nu),
\eeq
i.e., it defines an element of the first cohomology group of $\Lambda_0$ with
coefficients in the sheaf of holomorphic functions,
$f\in H^1_{gr}(\Lambda_0,H^0(\mathbb C^d, \O))$. Once again using the same
arguments, as that used in the proof of the part (b) of the Lemma 12 in \cite{shiota},
we get that there exists a holomorphic function $h_{s+1}(z,\nu)$
such that
\beq\label{new4}
f_{s+1}^{\l}(z,\nu)=h_{s+1}(z+\l,\nu)-h_{s+1}(z,\nu)+\wt B_{s+1,\,\nu}^{\l},
\eeq
where $\wt B_{s+1,\,\nu}^{\l,}$ is $z$-independent.
Hence, the function $\zeta_{s+1}=\chi_{s+1}+h_{s+1}$
has the following monodromy properties
\beq\label{new5}
\zeta_{s+1}(z+\l,\nu)-\zeta_{s+1}(z,\nu)=
\wt B_{s+1,\nu}^{\l}
+\sum_{i=1}^{s} B^{\,\lambda}_{i,\,\nu-1-s+i}\,\xi_{s+1-i}(z,\nu)\,.
\eeq
Let us consider the function $R_{s+1}$
defined by the equation
\beq\label{R}
R_{s+1}=\zeta_{s+1}(z-W,\nu)-\zeta_{s+1}(z+W,\nu)-
u(z,\nu)\,\xi_{s}(z,\nu-1)\,.
\eeq
Equation (\ref{may2}) and (\ref{may3}) imply that the r.h.s of (\ref{R}) has
no pole at $\TD\pm W$. Hence,
$R_{s+1}(z,\nu)$ is a holomorphic function of $z$.
From (\ref{new1},\ref{new5}) it follows that it is periodic
with respect to the lattice $\Lambda_W$, i.e., it is a function on $Y$.
Therefore, $R_{s+1}$ is a constant ($z$-independent) on each of the connected
components of $\C$.

Hence, the function
\beq\label{sol}
\xi_{s+1}(z,\nu)=\zeta_{s+1}(z,\nu)+l_{s+1}(z,\nu)+c_{s+1}(\nu)\,,
\eeq
where $c_{s+1}(\nu)$ is a constant, and $l_{s+1}$ is a linear form such that
$$l_{s+1}(2W,\nu)=-R_{s+1}(\nu)\,,$$
is a solution of (\ref{laxm2}). It satisfies the monodromy relations (\ref{new1})
with
\beq\label{new7}
B_{s+1,\,\nu}^{\,\l}=\wt B_{s+1,\,\nu}^{\,\l}+l_{s+1}(\l,\nu)\,.
\eeq
The induction step is completed and thus the lemma is proven.

On each step the ambiguity in the construction of $\xi_{s+1}$ is
a choice of linear form $l_{s+1}(z,\nu)$ and constants $c_{s+1}(\nu)$. As
in Section 2,
we would like to fix this ambiguity by normalizing monodromy coefficients
$B_{i,\,\nu}^{\l}$ for a set of linear independent vectors
$\l_0,\l_1,\ldots,\l_d\in\Lambda_W$.
It turns out that in the fully discrete case there is an obstruction for that.
This obstruction is a possibility of the existence of periodic solutions
of (\ref{laxm2}), $\xi_{s+1}(z+\l,\nu)=\xi_{s+1}(z,\nu), \ \l\in \Lambda_W,$
for $0\leq s\leq r-1$.

Note, that there are no periodic solutions of (\ref{laxm2}) for all $s$. Indeed,
the functions $\xi_s(z,\nu)$ as solutions of non-homogeneous equations are linear
independent. Suppose not.  Take a smallest nontrivial linear relation among
$\xi_s(z,\nu)$, and apply (5.24) to obtain a smaller linear relation.
The space of meromorphic functions on $Y$
with simple pole at $\TD$ is finite-dimensional. Hence, there exists minimal
$r$ such that equation (\ref{laxm2}) for $s=r$ has no periodic solutions.
\begin{lem} Let $\l_0,\l_1,\ldots,\l_d$ be a set of linear independent vectors
in $\Lambda_W$. Suppose equations (\ref{laxm2}) has periodic solutions for $s<r$
and has a quasi-periodic solution $\xi_{r}$ whose monodromy relations for
$\l_j$ have the form
\beq\label{new10}
\xi_{r}(z+\l_j,\nu)-\xi_{r}(z,\nu)=b^{\l_j},\ \ j=0,\ldots,d,
\eeq
where $b^{\l_i}$ are constants such that there is no linear form $l(z)$ on $Y$
with $l(\l_j)=b^{\l_j}$ and $l(2W)=0$.
Then for all $s$ equations (\ref{laxm2}) has solutions of the form (\ref{new})
satisfying (\ref{new1}) with $B_{i,\,\nu}^{\l_j}=b^{\l_j}\delta_{i,r}$,
i.e.,
\beq\label{new20}
\xi_{s}(z+\l_j,\nu)-\xi_{s}(z,\nu)=b^{\l_j}\xi_{s-r}(z,\nu).
\eeq
\end{lem}
{\it Proof.} We will now prove the lemma by induction in $s\geq r$. Let us assume
inductively that $\xi_{s-r}$ is known, and for $1\leq i\leq r$ there are solutions
$\tilde \xi_{s-r+i}$ of (\ref{laxm2}) satisfying
(\ref{new1}) with $B_{i,\,\nu}^{\l_j}=b^{\l_j}\delta_{i,r}$. Then,
according to the previous lemma, there exists a solution $\tilde\xi_{s+1}$
of (\ref{laxm2}) having the form (\ref{new}) and  satisfying monodromy relations
(\ref{new1}), which for $\l_j$ have the form
\beq\label{new30}
\tilde \xi_{s+1}(z+\lambda_j,\nu)-\tilde \xi_{s+1}(z,\nu)=
b^{\l_j}\tilde \xi_{s-r+1}(z,\nu)+
B_{s+1,\,\nu}^{\l_j}\,.
\eeq
If $\xi_{s-r}$ is fixed, then the general quasi-periodic solution $\xi_{s-r+1}$
with the normalized monodromy relations is of the form
\beq\label{transc}
\xi_{s-r+1}(z,\nu)=\wt \xi_{s-r+1}(z,\nu)+c_{s-r+1}(\nu)\,,
\eeq
where $c_{s-r+1}$ are constants on each component of $\C$.
It is easy to see that under the transformation (\ref{transc}) the functions
$\wt \xi_{s-r+i}$ get transformed to
\beq\label{new40}
\xi_{s-r+i}(z,\nu)=
\wt \xi_{s-r+i}(z,\nu)+c_{s-r+1}(\nu-i+1)\, \xi_{i-1}(z,\nu)\,.
\eeq
This transformation does not change the monodromy properties of $\xi_{s-r+i}$
for $i\leq r$, but changes the monodromy property of $\xi_{s+1}$:
\beq\label{new31}
\xi_{s+1}(z+\lambda_j,\nu)-\xi_{s+1}(z,\nu)=b^{\l_j}\xi_{s-r+1}(z,\nu)+
b^{\l_j}(c_{s-r+1}(\nu-r)-c_{s-r+1}(\nu))+B_{s+1,\,\nu}^{\l_j}
\eeq
Recall, that $\wt \xi_{s+1}$ was defined up to a linear form $l_{s+1}(z,\nu)$
which vanishes on $2W$. Therefore the normalization of
the monodromy relations for $\xi_{s+1}$ uniquely defines this form and
the differences $(c_{s-r+1}(\nu-r)-c_{s-r+1}(\nu))$. The induction
step is completed and the lemma is thus proven.

Note, the following important fact: if $\xi_{s-r}$ is fixed
then $\xi_{s-r+1}$, such that there exists quasi-periodic solution $\xi_{s+1}$
with normalized monodromy properties,
is defined uniquely up to the transformation:
\beq\label{transcc}
\xi_{s-r+1}(z,\nu)\longmapsto\xi_{s-r+1}(z,\nu)+c_{s-r+1}(\nu),
\ \ c_{s-r+1}(\nu+r)=c_{s-r+1}(\nu).
\eeq
Our next goal is to show that the assumption of Lemma 5.4 holds for some $r$, and
then to prove that the singular locus $\Sigma$ is in fact empty.

Shifting $z\to Z+z$, we get, as a direct corollary of Lemma 5.3,  that:\
if $Z\notin \cup_{i=0}^{s-1} (\Sigma_0-iV)$, where
$\Sigma_0=\bigcap_{\,k\in\mathbb Z}T_{U-V}^k\Theta\,, $
then there exist holomorphic functions $\tau_s(Z+z)$, which are {\it local} functions
of the variable $Z\in \mathbb C^g$ and {\it global} function of the variable
$z\in \C$, such that the equations
\beq\label{new50}
{\tau_{s}(Z)\over\theta(Z)}-{\tau_{s}(Z+U-V)\over\theta(Z+U-V)}=
{\theta(Z+U)\, \tau_{s-1}(Z-V)\over \theta(Z)\, \theta(Z+U-V)}, \
\eeq
holds, and the functions $\xi_s=\tau_s/\theta$ satisfy the monodromy relations
\beq\label{new60}
\xi_s(Z+z+\lambda)-\xi_{s}(Z+z)=
\sum_{i=1}^s B^{\,\lambda}_{i}(Z)\,\xi_{s-i}(Z+z)\,,\ \
\ \lambda\in \Lambda_W.
\eeq
If $\xi_{s-1}$ is fixed then $\xi_s$ is unique up to the transformations
\beq\label{new61}
\xi_s(Z+z)\to \xi_s(Z+z)+l_s(Z,z)+c_s(Z),
\eeq
where $l_s$ is a linear form in $z$ such that $l_s(Z,U-V)=0$, and $c_s(Z)$ are
$z$-independent.

Let $r$ be the minimal integer such that
$\xi_1,\ldots,\xi_{r-1}$ are periodic functions of $z$ with respect to
$\Lambda_W$, and there is no periodic solution $\xi_r$ of (\ref{new50}).
As it was noted above, the functions $\tau_s$ are linear independent. Hence,
$r\leq h^0(Y,\theta|_Y)$.

If $\xi_{r-1}$ is periodic, then the monodromy relation for
$\xi_r$ has the form
\beq\label{new65}
\xi_r(Z+z+\lambda)-\xi_{r}(Z+z)=B^{\,\lambda}_{r}(Z)\,,\ \
\ \lambda\in \Lambda_W.
\eeq
The function $B_r^\l$ is independent of the ambiguities in the definition
of $\xi_i, \ i<r$, and therefore, it
is a well-defined holomorphic function of $Z\in X$ outside of
the set $\cup_{i=1}^{r-1} (\Sigma-iV)$. The later is of codimension at least 2.
Hence, by Hartogs' theorem $B_{r}^\l(Z)$ extends to a holomorphic function
on $X$. Hence, it is a constant $B_{r}^\l(Z)=b^\l$. It was supposed that
the function $\xi_r$ can not be made periodic by transformation (\ref{new61}). Therefore,
there is no linear form on $\C$ such that $l(\l)=b^\l, \ l(U-V)=0$,
and the initial assumption of lemma 5.4 is proved.

\begin{lem} If equation (\ref{cm7d}) is satisfied, then
the singular locus $\Sigma\in \Theta$  is empty.
\end{lem}
{\it Proof.} The functions $\tau_1(Z+z)$ are defined as solutions of (\ref{new50}) along
$\C$. The restriction of $\tau_1(Z)$ on $\Theta$ is given by the formulae
(\ref{may2},\ \ref{may3}) for $s=0$, i.e.
\beq\label{may3a}
\tau_1={\theta(Z+U)\, \theta(Z-V)\over \theta(Z+U-V)}=
-{\theta(Z-U)\, \theta(Z+V)\over \theta(Z-U+V)}
\eeq

Let us first show that $\Sigma$ is invariant under the shift
by $rV$ (or equivalently by $rU$), where $r$ is defined above
(minimal integer such that there is no periodic solution $\xi_r$ of (\ref{new50}).
The functions $\tau_1(Z+z)$ are defined as solutions of (\ref{new50}) along
$\C$, and a priori there are no relations between $\tau_1$ defined for
$Z$ and its translates $Z-iV$. As shown in Lemma 5.4, the requirement that
there exists $\xi_{r+1}$ with normalized monodromy relations, defines $\tau_1$
uniquely, up to the transformations (\ref{new61}) with $l_1=0$ and with $rV$-periodic
$c_1$, i.e., $c_1(Z)=c_1(Z+rU)=c_1(Z+rV)$.

Let $Z$ be in $\Sigma$ and $Z+rV$ is not. Then $\tau_1$ can be defined as
a holomorphic function in the whole neighborhood of $(Z+rV)$. Therefore,
$\tau_1(Z)$ can be defined as a single-valued holomorphic function of $Z$ outside
of $\Sigma$. Hence, by Hartogs' arguments it can be extended across $\Sigma$.
The contradiction proves that $\Sigma=\Sigma+rV$.

By definition, $\Sigma$ is not invariant under the shift by $V$. Hence, it is
empty or $r>1$. Let $Z\in \Sigma$, then the r.h.s of equation (\ref{new50})
for $\tau_1(Z+V)$ vanishes. Therefore, $\xi_1$ is a constant along $\Sigma+V$.
Using the transformation (\ref{new61}) we can make it to be equal zero on $\Sigma+V$,
$\tau_1(Z+V)=0, \ Z\in \Sigma$. The same arguments applied consecutively
show that we may assume that $\tau_i(Z+iV)=0,\ \ i\leq r-2$. For $i=r-1$,
using in addition the equation $\theta(Z+rU)=0$ (which is due to the fact
$\Sigma=\Sigma+rU$), we get that, up to the transformation (\ref{new61}),
the function $\xi_{r-1}$ has vanishing order on $\Sigma+(r-1)V$ such
that the r.h.s of equation for $\xi_r$ on $\Sigma+rV$ is zero. Hence, $\xi_r$
can be defined as holomorphic function in the neighborhood of $\Sigma+rV$, and
restricted on $\Sigma+rV$ is a constant. That contradicts to the assumption
$b^\l\neq 0$, and thus the lemma is proven.

As shown above, if $\Sigma$ is empty, then the functions $\tau_s$ can be defined
as global holomorphic functions of $Z\in \mathbb C^g$. Then, as a corollary
of the previous lemmas we get the following statement.

\begin{lem} Let equation (\ref{cm7d}) for $\theta(Z)$ holds.
Then  there exists a formal solution
\beq\label{ff1}
\phi=1+\sum_{s=1}^{\infty}\xi_s(Z)\,k^{-s}
\eeq
of the equation
\beq\label{ff2}
k\phi(Z+V,k)=k\phi(Z+U, k)+u(Z)\,\phi(Z, k)\,,
\eeq
with
\beq\label{ff3}
u(Z)={\theta(Z+U+V)\,\theta(Z)\over \theta(Z+U)\,\theta(Z+V)},
\eeq
such that:

(i) the coefficients $\xi_s$ of the formal series $\phi$
are of the form $\xi_s=\tau_s/\theta$, where
$\tau_s(Z)$ are holomorphic functions;

(ii) $\phi(Z,k)$ is quasi-periodic with respect to the lattice $\Lambda$ and
for the basis vectors $\l_j$ in $C^g$ its monodromy relations
have the form
\beq\label{ff4}
\phi(Z+\lambda_j)=(1+b^{\l_j}\,k^{-1})\,\phi(Z,k),\ \ j=1,\ldots, g,
\eeq
where $b^{\l_j}$ are constants such that there is no linear form on $\mathbb C^g$
vanishing at $\l_j$ and $U-V$, i.e., $l(\l_j)=\l(U-V)=0$;

(iii) $\phi$ is unique up to the multiplication by a constant in $Z$ factor.
\end{lem}

\bigskip
\noindent
{\bf Commuting difference operators.} The formal series $\phi(Z,k)$ defines a unique
pseudo-difference operator
\beq\label{LLd}
\L(Z)=T+\sum_{s=0}^{\infty} w_s(Z)\,T^{-s},\ \ T=e^{\p_m},
\eeq
such that the equation
\beq\label{kkd}
\left(T+\sum_{s=0}^{N} w_s(Z+mU+nV)\,T^{-s}\right)\,\psi=k\psi\,.
\eeq
holds. Here $\psi=k^{n+m}\phi(nV+mU+Z,k)$.
The coefficients \, $w_s(Z)$ of $\L$  are meromorphic functions on the
abelian variety $X$ with poles along  the divisors $T_U^{-i}\Theta= \Theta-iU,\ i\leq s+1$.

From equations (\ref{ff2},\,\ref{kkd}) it follows that
\beq\label{m22}
\left((\Delta_1\L^i)\,T_1-(\Delta \L^i) \,T-[u,\L^i]\right)\psi=0
\eeq
where $\Delta_1\L^i$ and $\Delta\L^i$ are pseudo-difference operator in $T$, whose
coefficients are difference derivatives of the coefficients of $\L^i$
in the variables $n$ and $m$ respectively.
Using the equation $(T_1-T-u)\,\psi=0$, we get
\beq\label{m23}
\left(\left(\Delta_1\L^i\right)T-\left(\Delta\L^i\right) T +
\left(\Delta_1 \L^i\right)u-[u,\L^i]\right)\psi=0.
\eeq
The operator in the left hand side of (\ref{m23}) is a pseudo-difference operator
in the variable $m$. Therefore, it has to be equal to zero. Hence, we have the equation
\beq\label{m21}
\left(\Delta_0\L^i\right) T +\left(\Delta_V \L^i\right)u-[u,\L^i]=0, \
\Delta_{0}=e^{\p_V}-e^{\p_U}
\eeq
As before, the strictly positive difference part of the operator $\L^i$ we denote by
$\L^i_+$. Then,
\beq\label{m24}
\left(\Delta_{0}\L^i_+\right) T +\left(\Delta_V \L^i_+\right)u-[u,\L^i_+]=
-\left(\Delta_{0}\L^i_-\right) T -\left(\Delta_V \L^i_-\right)u+[u,\L^i_-]
\eeq
The left hand side of (\ref{m24})
is a difference operator with non-vanishing coefficients only at the positive
powers of $T$. The right hand side is a pseudo-difference operator
of order $1$. Therefore, it has the form $f_i T$. The coefficient $f_i$ is
easy expressed in terms of the leading coefficient $\L^i_-$. Finally we get
the equation
\beq\label{m25}
\left(\Delta_{0}\L^i_+\right) T +\left(\Delta_V \L^i_+\right)u-[u,\L^i_+]=
-(\Delta_{0}F_i)\,T,
\eeq
where $F_i=\res\  \L^i$. The vanishing
of the coefficient at $T^0$ in the right hand side of (\ref{m24}) implies the equation
\beq\label{lax55d}
\Delta_{0} F_i^1=-\left(\Delta_V F_i\right) u,\ \ F_i^1=\res \ \L^iT\,,
\eeq
analogous to (\ref{lax55}).

\begin{lem} The abelian functions $F_i$ have at most simple poles on the divisors
$\Theta$ and $\Theta_U$.
\end{lem}
The wave solution $\psi$ define the unique operator $\Phi$ such that
\beq\label{Sd}
\psi=\Phi k^{n+m},\ \
\Phi=1+\sum_{s=1}^{\infty}\f_s(Um+Vn+Z)\,T^{-s}.
\eeq
The dual wave function
\beq\label{psinewd}
\psi^+=k^{-n-m}\left(1+\sum_{s=1}^{\infty}\xi^+_s(Um+Vn+Z)\, k^{-s}\right)
\eeq
is defined by the formula
\beq\label{m26}
\psi^+=k^{-n-m}\, T_1\,\Phi^{-1}\,T_1^{-1}.
\eeq
It satisfies the equation
\beq\label{adjd}
(T_1^{-1}-T^{-1}-u)\,\psi^+=0,
\eeq
which implies that the functions $\xi_s^+(Z)$ have
the form $\xi_s^+(Z)=\tau_s^+(Z)/\theta(Z+U+V)$,
where $\tau_s^+$ are holomorphic functions.
Therefore, the functions $J_s(Z)$ such that
\beq\label{Jnew}
(\psi^+T_1)\,\psi=k+\sum_{s=1}^{\infty}J_s(Um+Vn+Z)\,k^{-s+1}
\eeq
are meromorphic function on $X$ with the simple poles at $\Theta$ and
$T_U^{-1}\Theta=\Theta_U$.

From the definition of $\L$ it follows that
\beq\label{20new}
\res_k\left((\psi^+T_1)\,(\L^n\psi)\right)k^{-2}dk=
\res_k\left((\psi^+\,T_1)\,\psi\right)k^{n-2}dk
=J_{n}.
\eeq
On the other hand, using the identity (\ref{dic}) we get
\beq\label{20nn}
\res_k((\psi^+\,T_1)\,\L^n\psi)\,k^{-2}dk=
\res_k\left(k^{-n-m}\Phi^{-1}\right)\left(\L^n\Phi k^{n+m}\right)k^{-1}dk=
\res_{T}\L^n=F_n.
\eeq
Therefore, $F_n=J_{n}$ and the lemma is proved.

The rest of the proof of Theorem 1.2 is identical to that in the proof
of Theorem 1.1. Namely: lemma 5.7 directly implies that for the generic $Z\in X$
linear combinations of operators $\L^i_+$
span commutative rings $\A^Z$ of ordinary difference operators.
They define a spectral curve $\G$ with two smooth points $P_{\pm}$ and a map (\ref{is}).
The global existence of the wave function implies
equations (\ref{nnov7}).

Equations (\ref{z4}) define the KP hierarchy deformations of these rings.
From (\ref{z4},\,\ref{kkk}), and equation (\ref{laxm2}) for $s=0$ we get
\beq\label{m30}
\p_{t_n}w_0=\Delta_U F_n^1,\ \ w_0=-\Delta_U\xi_1,\ u=\Delta_0 \xi_1\,.
\eeq
Then, (\ref{nnov7}) and (\ref{lax55d}) imply
\beq\label{m31}
\p_{\,t_n}\ln u=\Delta_U\Delta_V (\p_{V_n}\ln \theta).
\eeq
By definition, $u$ is given by the formula (\ref{ud}), i.e.,
$\ln u=\Delta_U\Delta_V\ln \theta$. Therefore, equation (\ref{m31})
identifies $\p_{t_n}$ with $\p_{V_n}$. Hence,
the orbit of the KP flows is in $X$. Hence the generalized Jacobian $J(\G)$
of the spectral curve is compact and the spectral curve is smooth. As in the previous case
these arguments complete the proof of the theorem.

\end{document}